\documentclass[a4paper,11pt]{article}

\title{A nonsmooth nonconvex descent algorithm}

\author{ Jan Mankau, Friedemann Schuricht \\[2mm]
         {\fns TU Dresden, Fakult\"at Mathematik} \\[-1mm] 
         {\fns 01062 Dresden, Germany} 
       }

\usepackage[a4paper,left=26mm,top=25mm,right=26mm,bottom=34mm]{geometry}
\date{}
\usepackage{amssymb}
\usepackage{amsmath}
\usepackage{amsthm}
\usepackage{mathrsfs}  

\usepackage{marginnote}                
           \reversemarginpar           

\usepackage{ifthen}
\usepackage{eurosym}

\usepackage{varioref}        

\usepackage{pdfpages}

\numberwithin{equation}{section}      

\newtheorem{theo}[equation]{Theorem} 
\newtheorem{prop}[equation]{Proposition}

\newtheorem{lem}[equation]{Lemma}
\newtheorem{cor}[equation]{Corollary}

\theoremstyle{definition} 

\newtheorem{rem}[equation]{Remark}  
\newtheorem{exam}[equation]{Example}  
\newtheorem{ass}[equation]{Assumption}
\newtheorem{alg}[equation]{Algorithm}

\newcounter{nr}  \labelformat{nr}{{\normalfont \arabic{nr}}}

\newcommand{\bgl}[1][{\normalfont (\arabic{nr})}]
   {\begin{list} {#1} {\usecounter{nr}
          \setlength{\topsep}{0.5ex plus0.2ex minus0.1ex}     
          \setlength{\itemsep}{0.2ex plus0.05ex minus0.03ex}  
          \parsep0pt \itemindent0pt 
          \leftmargin30pt   \labelwidth20pt}   }
\newcommand{\el}{\end{list}}

\newcounter{nra}  \labelformat{nra}{{\normalfont \alph{nr}}}
\newcommand{\bgla}[1][{\normalfont (\alph{nra})}]
   {\begin{list} {#1} {\usecounter{nra}
          \setlength{\topsep}{0.5ex plus0.2ex minus0.1ex}     
          \setlength{\itemsep}{0.2ex plus0.05ex minus0.03ex}  
          \parsep0pt \itemindent0pt 
          \leftmargin30pt   \labelwidth20pt}   }
\newcommand{\ela}{\end{list}}

\newcommand{\bgln}                
  {\bg{list}{(\arabic{section}.\arabic{equation})}{\usecounter{equation} 
          \setlength{\topsep}{0.5ex plus0.2ex minus0.1ex}     
          \setlength{\itemsep}{0.2ex plus0.05ex minus0.03ex}  
          \parsep0pt \itemindent0pt 
          \leftmargin30pt   \labelwidth20pt}   }

\newcommand{\beq}{\begin{equation}\nonumber}
\newcommand{\bee}{\begin{equation}\nonumber}
\newcommand{\bn}[1]{\begin{equation}\label{#1}} 
\newcommand{\bnn}[1]{\begin{equation}\llb{#1}\nonumber}

\newcommand{\ee}{\end{equation}}

\newcommand{\ba}{\begin{eqnarray}} \newcommand{\ea}{\end{eqnarray}}
\newcommand{\barr}{\begin{array}} \newcommand{\earr}{\end{array}}

\newcommand{\nn}{\nonumber}

\newcommand{\ben}[1][]{\begin{equation}\label{#1}}
\newcommand{\benn}[1][]{\begin{equation}\llb{#1}\nonumber}

\newcommand{\en}{\end{equation}}

\newenvironment{pf}[1]{\noindent {\sc Proof}#1.\hspace{1pt}}%
     {\hspace*{\fill} $\diamondsuit$\par\medskip}

\newcommand{\bg}[1]{\begin{#1}} \newcommand{\e}[1]{\end{#1}}

\newcommand{\z}{\enspace}    

\newcommand{\hs}[1]{\leer\hspace{#1ex}}
\newcommand{\qmz}[1]{\quad\mbox{#1}\z}
\newcommand{\qmq}[1]{\quad\mbox{#1}\quad}
\newcommand{\zmz}[1]{\z\mbox{#1}\z}
\newcommand{\leer}{\mbox{}}
\newcommand{\tx}[1]{\text{#1}}

\newcommand{\msk}{\medskip}
\newcommand{\bsk}{\bigskip}   

\newcommand{\fns}{\footnotesize}

\newcommand{\noi}{\noindent}
  \newcommand{\ol}{\overline}

\newcommand{\reff}[1]{{\rm (\ref{#1})}}

\newcommand{\bmp}[1]{\begin{minipage}[t]{#1}}
\newcommand{\bmpc}[1]{\begin{minipage}[c]{#1}}
\newcommand{\emp}{\end{minipage}}

\newcommand{\op}[1]{\operatorname{#1}}   

\newcommand{\osm}[2]{\ensuremath{\overset{#1}{#2}}}

\newcommand{\wto}{\rightharpoonup}

\newcommand{\lto}{\ensuremath{\longrightarrow}}

\newcommand{\dl}{\langle}  \newcommand{\dr}{\rangle}

\newcommand{\pa}{\partial}

\newcommand{\subs}{\subset}

\renewcommand{\d}[1][]{{\rm d}}     

\newcommand{\ti}{\tilde}

\newcommand{\N}{\ensuremath{\mathbb{N}}} 
 
\newcommand{\R}{\ensuremath{\mathbb{R}}}

\def\XXint#1#2#3{{\setbox0=\hbox{$#1{#2#3}{\int}$}
\vcenter{\hbox{$#2#3$}}\kern-.5\wd0}}

\newcommand{\resp}{respectively}

\newcommand{\dst}{\displaystyle}

\newcommand{\al}{\alpha}

\newcommand{\ga}{\gamma}      
\newcommand{\de}{\delta}      
\newcommand{\ep}{\varepsilon}

\newcommand{\la}{\lambda}     \newcommand{\La}{\Lambda}

\newcommand{\si}{\sigma}      
\newcommand{\ta}{\tau}

\newcommand{\gG}[3]{\partial^{#3}\!#1(#2)}
\newcommand{\gdD}[4]{#1 ^0_{#2}( #3; #4) }

\newcommand{\DefFunktion}[3]{#1:\,#2\to #3 }

\newcommand{\tme}{\:\big|\:}    
\newcommand{\sk}[3]{\langle #1, #2\rangle_{#3}} 
\newcommand{\norm}[1]{\| #1 \| }       
\newcommand{\norms}[2]{\left\| #1 \right\|_{#2} }       
\newcommand{\set}[2]{\left\{ #1 \tme #2\right\} }               
\newcommand{\sset}[1]{\left\{ #1 \right\} }               
\newcommand{\ball}[3]{ B_{#2}(#1) }         

\newcommand{\cball}[3]{\overline {B_{#2}(#1)}} 

\newcommand{\USg}{G}
\newcommand{\ThN}{T_1}
\newcommand{\NHI}{T_2}
\newcommand{\fkt}{T_3}

\usepackage{graphicx}
\usepackage{multirow}
\usepackage{enumerate}   
\usepackage{float}

\allowdisplaybreaks

\begin{document}

\bibliographystyle{plain}


\maketitle


\bg{abstract} The paper presents a new descent algorithm for locally Lipschitz
continuous functions $f:X\to\R$.  The selection of a descent direction at some
iteration point $x$ combines an approximation of the set-valued gradient of
$f$ on a suitable neighborhood of $x$ (recently introduced by 
Mankau \& Schuricht
\cite{01-paper}) with an Armijo type step control. The algorithm is 
analytically justified and it is shown that accumulation points of iteration
points are critical points of $f$. Finally the algorithm is tested for
numerous benchmark problems and the results are compared with simulations
found in the literature. 
\e{abstract}

\section{Introduction}
In this paper we present a new descent algorithm to find local minima or
critical points of a locally Lipschitz continuous function $f:X\to\R$ on a
Hilbert space $X$. For the minimization of a nonsmooth function  
\begin{equation}\label{P} 
 f(x)\to\min\;\;\;\;(x\in X)
\end{equation}
numerous algorithms based on quite different methods have been proposed in the
literature. Let us mention, without being complete,  
bundle-type methods 
(cf. Alt \cite{Alt2004}, Frangioni \cite{Frangioni2002}, Gaudioso \& Monaco
\cite{Gaudioso1982}, Hiriat-Urruty \cite{Hiriat1993}, Kiwiel \cite{Kiwiel1985},
Makela \& Neittaanmaki \cite{Makela1992}, Mifflin \cite{Mifflin1977},
Schramm \cite{Schramm1989}, Wolfe \cite{Wolfe1975}, Zowe \cite{Zowe1985}),
proximal point and splitting methods as e.g. the Fista or the primal dual
method (cf. Beck \cite{Beck2017}, Eckstein \& Bertsekas \cite{EckBert1992},
Chambolle \& Pock \cite{ChaPock2011}),
gradient sampling algorithms (cf. 
Burke, Lewis \& Overton \cite{CS-Verfahren-Burke},
Kiwiel \cite{CS-Verfahren-Kiwiel}),
algorithms based on smoothing techniques (cf. Polak \& Royset
\cite{Polak2003}) and the discrete gradient method 
(cf. Bagirov \& Karasozen \cite{Bagirov2008}).

Bundle-type methods, proximal point methods, and splitting methods require 
$f$ to be convex or to have some other special structure. 
Many algorithms for locally Lipschitz continuous functions 
as the discrete gradient method need to know the entire generalized gradient
of $f$ at given points.  
Stochastic methods like the gradient sampling algorithm are robust without 
the knowledge of the entire generalized gradient, but at the cost of high
computational effort. Therefore they are limited to minimization problems on low
dimensional spaces. 

Recall that the derivative $f'(x)$ indicates a direction of descent for $f$
near $x$. However if the direction of descent changes rapidly in a small
neighborhood of $x$, which is typical for functions $f$ having large second
derivatives or that are even nonsmooth, then some knowledge of $f$ on a whole
neighborhood of $x$ 
is necessary for the determination of a suitable direction of
descent near $x$.  

For a new robust and fast algorithm we combine ideas of bundle-methods and 
gradient sampling methods. We use the concept of gradients of $f$ 
on sets as introduced in Mankau \& Schuricht~\cite{01-paper} 
(which extends ideas from Goldstein \cite{Goldstein}).
Here, similar to gradient sampling methods, generalized gradients of $f$ on a
whole neighborhood of $x$ are considered for the determination of a
suitable descent direction near point $x$. But, in contrast to e.g. 
Burke, Lewis \& Overton \cite{CS-Verfahren-Burke} and Kiwiel
\cite{CS-Verfahren-Kiwiel}, the set-valued gradient 
on a neighborhood of point $x$ is not approximated stochastically. 
We rather use an elaborate recursive inner approximation coupled with the 
computation of related descent directions until a generalized Armijo condition
is satisfied (a condition similar to that used in Alt \cite{Alt2004} and 
Schramm \cite{Schramm1989} in connection with the 
$\ep$-subdifferential). Finally a line search along a direction of
sufficient descent gives the next iteration point 
(cf. Pytlak \cite{Pyt2009}). 
For better performance we may also adapt the norm of $X$ in each step.
It turns out that our algorithm requires  
substantially less gradient computations 
than in \cite{CS-Verfahren-Burke} and \cite{CS-Verfahren-Kiwiel}.
Therefore it is also applicable on high dimensional spaces as needed for
variational problems. 

For a locally Lipschitz continuous function $f:X\to\R$ our methods merely
demand that, at any point $x$, both the value $f(x)$ and at least  
one element of the generalized gradient $\partial f(x)$ (in the sense of
Clarke \cite{Clarke}) can be computed. Notice that this mild requirement 
is assumed in any of the above mentioned gradient based algorithms and that it
is typically met in applications. In an upcoming paper an extended 
algorithm is presented where quasi-Newton methods and preconditioning methods 
are included by a suitable change of norm in each iteration step. 

Section~\ref{gg} gives a brief overview about gradients on sets as needed for
our treatment. The algorithm and some convergence
results are given in Section~\ref{SecAlg}. After the formulation of 
the condition 
of sufficient descent and of several general
assumptions, 
Section~\ref{SecAlg1} provides the main Algorithm \ref{A-Main}
and its properties. Algorithm~\ref{A-Main} calls 
Algorithm~\ref{A-IA} for the computation of a suitable inner approximation of
the set-valued gradient on a neighborhood of the current iteration point 
and the computation of a related descent direction while Step 3 of
Algorithm~\ref{A-IA} calls 
Algorithm~\ref{A-IS} for some subiteration. Figure \ref{MAlg-Flussdiagramm} 
gives an overview of the whole algorithm and several statements
justify essential steps of it. Theorem~\ref{T-Main}
shows that every accumulation point of iteration points produced by
Algorithm~\ref{A-Main} is a critical point in the sense of Clarke. 
The proofs are collected in Section~\ref{Proofs}.
Comprehensive numerical tests of our algorithm for classical
benchmark problems can be found in Section~\ref{Benchmark}. Here the 
simulations are also compared with results from 
Burke, Lewis \& Overton \cite{CS-Verfahren-Burke}, 
Kiwiel \cite{CS-Verfahren-Kiwiel}, Alt \cite{Alt2004},  
Schramm \cite{Schramm1989} and the BFGS algorithm. 

\msk

\noi
{\it Notation:} $X$ is a Hilbert 
space\footnote{Notice that any Hilbert space is uniformly convex and
reflexive.} 
with scalar product $\dl\cdot,\cdot\dr$ where the dual $X^*$ is always 
identified with $X$.
For a set $M$ we write $\ol{M}$ for its closure, 
$\op{conv}M$ for its convex hull and 
$\overline{\op{conv}}\, M$ for its closed convex hull.
$B_\ep(x)$ and $B_\ep(M)$ are the open $\ep$-neighborhood of point $x$ and
set $M$, \resp. $]x,y[$ stands for the open segment between $x,y$ and, in
particular, $]x,y[\,\subs\R$ for the open interval.
$\R_{>0}$ denotes the positive real numbers.
For a locally Lipschitz continuous function $f:X\to\R$ we write 
$f^0(x;h)$ for Clarke's generalized directional derivative of $f$ at $x$ in
direction $h$ and $\pa f(x)\subs X$ for Clarke's generalized gradient of $f$ at
$x$ (cf. Clarke~\cite{Clarke}).

\section{Gradient on sets}
\label{gg}

Let $\DefFunktion{f}{X}{\R}$ be a locally Lipschitz continuous function
on a Hilbert space $X$. Clarke's generalized gradient $\gG{f}{x}{}$
of $f$ at $x$ and the corresponding generalized directional derivative 
$\gdD{f}{}{x}{h}$ of $f$ at $x$ in direction $h$ 
somehow express the behavior of $f$ at point $x$ (cf. Clarke~\cite{Clarke}). 
However, for the construction of a descent step in a numerical scheme,  
some information about the behavior of $f$ on a whole
neighborhood of $x$ is useful in general. 
In particular, for describing the behavior of $f$  
on the whole $\ep$-ball $\cball{x}{\ep}{X}$,
we use some set-valued gradient $\gG{f}{x}{\ep}$ of $f$
and some corresponding generalized directional derivative 
$\gdD{f}{\ep}{x}{h}$ as
introduced in Mankau \& Schuricht~\cite{01-paper}
by using Clarke's pointwise quantities. For the convenience of the
reader we present some brief specialized summary of that material as 
needed for our treatment.   

For $\ep>0$ we define the {\it gradient} of $f$ on $\cball{x}{\ep}{X}$ by
\begin{equation}\label{gg-e1}
 \gG{f}{x}{\ep}:=\overline{\op{conv}}\bigcup\limits_{y\in\cball{x}{\ep}{X} } 
 \gG{f}{y}{}
\end{equation}
(notice that the closed convex hull $\ol{\op{conv}}$ agrees with the weakly 
closed convex hull) 
and the {\it directional derivative} of $f$ on $\cball{x}{\ep}{X}$
in direction $h\in X$ by 
\begin{equation}\label{gg-e2}
 \gdD{f}{\ep}{x}{h}:=\sup\limits_{y\in\cball{x}{\ep}{X} }\gdD{f}{}{y}{h}\,.
\end{equation}
We have the following basic properties (cf. 
Propositon~2.3 and Corollary~2.10 in \cite{01-paper}).

\begin{prop}\label{gg-s1}
Let $\DefFunktion{f}{X}{\R}$ be Lipschitz continuous of rank~$L$ on a
neighborhood of $\cball{x}{\ep}{X}$ with $x\in X$ and $\ep>0$. Then 
\bgl
 \item $\gG{f}{x}{\ep}$ is nonempty, convex, weak-compact and bounded by
   $L$. 
 \item $\gdD{f}{\ep}{x}{\cdot}$ is finite, positively homogeneous, subadditive,
   and Lipschitz continuous of rank $L$. Moreover it is the support function
   of $\gG{f}{x}{\ep}$, i.e. 
\ben
\gdD{f}{\ep}{x}{h} = \max\limits_{a\in\gG{f}{x}{\ep}} 
 \sk{a}{h}{}  \qmz{for all} h\in X\,. \label{supportgG}
\ee
 \item We have 
\ben
  \gG{f}{x}{\ep} = \set{a\in X}{\sk{a}{h}{}\leq  \gdD{f}{\ep}{x}{h}\zmz{for
      all} h\in X }\,.\label{ChargG}
\ee
\item \label{gg-s1-4}
  Let $h\in X$ with $\gdD{f}{\ep}{x}{h}<0$ and let $t>0$
  with $]x,x+th[\:\subset \cball{x}{\ep}{X}$. Then 
$$f(x+th)\leq f(x)+t\gdD{f}{\ep}{x}{h}<f(x)\; .$$
\item Let $\ep_k\to0$ with $\ep_k>0$ and let $h\in X$. Then 
\begin{equation*}
 \lim_{k\to\infty}\gdD{f}{\ep_k}{x}{h}=\gdD{f}{}{x}{h}
 \qmq{and} \bigcap_{k\in\N}\gG{f}{x}{\ep_k}=\gG{f}{x}{}\,.  
\end{equation*}
\el 
\end{prop}

\noindent

Regularity of $f$ at $x$, i.e.
$0\notin\gG{f}{x}{}$, implies regularity of $f$ on 
some $\cball{x}{\ep}{X}$ by Proposition~2.16 in \cite{01-paper}.
\begin{lem}\label{L-BW}
Let $f:X\to\R$ be locally Lipschitz continuous and let 
$0\notin\gG{f}{x}{}$ for some $x\in X$. 
Then there exist $\ep>0$ and $h\in X$ with $\norm{h}{}=1$ such that
\begin{equation*}
  -\norms{a}{}\leq \sk{a}{h}{}\leq \gdD{f}{\ep}{x}{h}<0
  \qmz{for all} a\in\gG{f}{x}{\ep}\,.
\end{equation*}
Moreover, $0\notin\gG{f}{x}{\ep}$ by {\rm (\ref{ChargG})}.
\end{lem}

Motivated by Proposition \ref{gg-s1} \reff{gg-s1-4}
we say that $h\in X$ is a {\it descent direction} of $f$ on $\cball{x}{\ep}{X}$
if $\gdD{f}{\ep}{x}{h}<0$.
We call $\ti h$ {\it steepest} or {\it optimal} descent
direction of $f$ on $\cball{x}{\ep}{X}$ (with respect to $\norm{\cdot}$) if
\begin{equation}
  \|\ti h\|=1 \qmq{and}
  \gdD{f}{\ep}{x}{\ti h}=\min\limits_{\norm{h}{}\leq1}  \gdD{f}{\ep}{x}{h}<0\,.
\end{equation}
Theorem 3.10 of \cite{01-paper} ensures the existence of optimal descent
directions and of norm-minimal elements in $\gG{f}{x}{\ep}$.

\begin{prop}\label{Satz1Umg}
Let $\DefFunktion{f}{X}{\R}$ be Lipschitz continuous on a neighborhood of 
$\cball{x}{\ep}{X}$ for some $x\in X$, $\ep>0$.
Then there is a unique $\ti a\in \gG{f}{x}{\ep}$ with 
\begin{equation}\label{deff_0}
 \norm{\ti a}{}=\min\limits_{a\in \gG{f}{x}{\ep}}\norm{a}{}\,.
\end{equation}
If $0\notin \gG{f}{x}{\ep}$ or, equivalently, 
$\gdD{f}{\ep}{x}{h}<0$ for some $h\in X$ (cf. \reff{ChargG}), then 
there is a unique optimal descent direction $\ti h$ on $\cball{x}{\ep}{X}$.
In particular
\begin{equation}\label{WeftverR}
  \ti h=-\frac{\ti a}{\norm{\ti a}{} }\,, \qquad 
  \gdD{f}{\ep}{x}{\ti h}=\min\limits_{ \norm{h}{}\leq1}  \gdD{f}{\ep}{x}{h}
  =-\norm{\ti a}{}\,.
\end{equation}
\end{prop}

\noindent
Corollary 3.15 and Corollary 3.16 in \cite{01-paper} state some stability of descent directions.
\begin{cor}
Let $\DefFunktion{f}{X}{\R}$ be Lipschitz continuous of rank $L$ on a
neighborhood of $\cball{x}{\ep}{X}$ for some $x\in X$, $\ep>0$, let 
$0\notin \gG{f}{x}{\ep}$, and let $\ti a$, $\ti h$ be as in 
Proposition~{\rm \ref{Satz1Umg}}. Then every $h\in X$ with 
$\norm{h-\ti h}{}<\frac{\norm{\ti a}{}}L$
is a descent direction on $\ball{x}{\ep}{X}$. 
\end{cor}

\noindent
This allows to get descent directions by suitable approximations of
$\ti a$, which is important for our numerical algorithms. 

\begin{cor}\label{Kiwiel-GS-Verfahren-Corollary}
Let $\DefFunktion{f}{X}{\R}$ be Lipschitz continuous on a
neighborhood of $\cball{x}{\ep}{X}$ for some $x\in X$, $\ep>0$, let 
$0\notin \gG{f}{x}{\ep}$, and let $\ti a$, $\ti h$ be as in 
Proposition~{\rm \ref{Satz1Umg}}.
Then for any $\de\in]0,1[$ there is some $\tau>0$
such that for every $a'\in\gG{f}{x}{\ep}$ with 
$$
  \norm{a'}{}\leq \min\limits_{  a\in \gG{f}{x}{\ep} }\norm{a}{}+\ta=
  \norm{\ti a}{}+\ta
$$ 
we have that $h':=\frac{-a'}{\norm{a'}{}}$ is a descent direction 
on $\cball{x}{\ep}{X}$ and satisfies
\begin{equation*}
  \gdD{f}{\ep}{x}{h'} \stackrel{(\ref{supportgG})}{=}\max\limits_{a\in
    \gG{f}{x}{\ep}} \sk{a}{h'}{}<- \de\norm{\ti a}\,.
\end{equation*}
\end{cor}

\section{Descent algorithm}\label{SecAlg}

We now introduce some descent algorithm for locally Lipschitz continuous
functions $f:X\to\R$ on a Hilbert space $X$. At each iteration point $x$
we determine an approximation $a$ of the norm-minimal element 
$\ti a\in\gG{f}{x}{\ep}$ (cf. (\ref{deff_0}))  
with respect to some suitable radius~$\ep>0$. We are interested in
pairs $(a,\ep)$ satisfying a {\it condition of sufficient descent}  
in the sense of a generalized Armijo step of the form
\begin{equation}\label{L-DSAs-0}
f(x-\ep h)-f(x)\leq -\de\ep\|a\| \qmq{with}
h=\frac{a}{\norms{a}{}} 
\end{equation}
where $\de\in\,]0,1[$ is fixed for the whole scheme. 
As new iteration point we then select $x-\si h$ for some $\si\ge\ep$  
such that \reff{L-DSAs-0} still holds with $\si$ instead of $\ep$.
If $0\in\gG{f}{x}{\ep}$, the norm $\|a\|$ will be very small and the 
{\it null step condition} 
\begin{equation}\label{L-NSAs-0}
 \|a\|<\ThN(\ep)
\end{equation}
(with a suitable control function $\ThN$ that is fixed for the whole scheme)  
indicates that situation. Here we cannot expect \reff{L-DSAs-0} in general 
and we have two possibilities. 
If $\ep$ is on the desired level of accuracy for the minimizer 
(or critical point), we can stop the algorithm. 
Otherwise the used ball $\cball{x}{\ep}{X}$ is too large 
for an iteration step with sufficient descent. Therefore we decrease
$\ep$ and look for sufficient descent with a new pair $(a,\ep)$. 
Our approximation of $\ti a$ combined with the analytically justified step size
control ensures that we always get sufficient descent for some $\ep>0$ 
small enough
(cf. Lemma~\ref{L-BW} and also the proof of Theorem~\ref{T-Main}).
That we finally end up with the null step condition on the desired scale,
$\ep$ has to become sufficiently small during the algorithm, which 
is ensured by control functions $\ThN$ and $\NHI$. But, that the algorithm
doesn't get stuck in a small ball without critical point, 
$\ep$ shouldn't approach zero too fast, which is ensured by
control functions $\ThN$ and $G$. Thus a careful selection of the step size,
that is related to $\ep$, plays a very important role.
The algorithm can be improved by choosing 
suitable equivalent norms at each iteration step. 

Let us start with general requirements for the control functions
$\ThN$, $\NHI$, and $G$.

\begin{ass}\label{Ass-1} Suppose that:
\bgl
\item\label{As-1a} $\ThN,\NHI:\R_{>0}\to\R_{>0}$ are non-decreasing functions
such that
\beq
\lim\limits_{t\to0}\ThN(t)=0 \qmq{and}
\lim\limits_{n\to\infty}\NHI^n(t)=0 \qmz{for all} t>0\,
\ee
where $T^n_2$ is given inductively by $T^1_2=T_2$ and $T^n_2=T_2\circ T^{n-1}_2$
for all $n\in\N$. 
Notice that this implies
\beq
T_2(t)<t \qmz{for all} t>0\,
\ee
(otherwise $\NHI(t)\geq t$ for some $t>0$ and, since 
$\NHI$ is non-decreasing, induction would give  
$\NHI^n(t)\geq t\stackrel{n\to\infty}{\nrightarrow}0$).

\item $\USg:\R_{>0}^2\to\R_{>0}$ is a function having at least one of the
following properties: 

\vspace {-2mm} 

\bgla
\item $\USg(x_k,y_k)\to 0$ implies $x_k\to 0$ for any sequences 
$(x_k)$ and $(y_k)$.

\item For any $x_0>0$ there is some $y_0$ such that 
$\USg(x,y)\geq y$ for all $x>x_0$ and $y<y_0$.
\ela

\el
\end{ass}

\noindent
Since the conditions for $T_2$ and $G$ are quite technical, we 
provide some typical examples.

\begin{exam}[examples for $T_2$ and $\USg$]$\;$    
\bgl
  \item $\NHI(x):=\al x$ for $\al\in\left]0,1\right[$\,.
  \item $\NHI(x):=\frac x{1+x}$ where, in particular, 
        $\NHI\big(\frac1n\big)=\frac 1{n+1}$\,.

  \item $\USg(x,y)\geq \al$ with a constant $\al>0$ satisfies (a).
  \item $\USg(x,y)=\al y$ with $\al\geq1$ satisfies (b).
  \item $\USg(x,y)=\overline{\USg}(x)$ with
        $\overline{\USg}:\R_{>0}\to\R_{>0}$ non-decreasing and
        $\lim\limits_{t\to0}\overline{\USg}(t)=0$ satisfies (a).
  \item $\USg:=\USg_1+\USg_2$ and $\USg:=\max\sset{\USg_1,\USg_2}$ satisfy 
        (a) or (b) if $\USg_1$, $\USg_2$ satisfy both (a) or (b),
        respectively. 
  \item $\USg(x,y):=\max\sset{\al,y}$ satisfies (a) and 
        $\USg(x,y):=\min\sset{\al,y}$ satisfies (b) for any $\al>0$.
\el
\end{exam}

As already mentioned, it might be useful to adapt the norm 
in every iteration 
(recall that the Newton method can be considered as descent algorithm 
with changing norm at each step). In our algorithm we allow a change of norm in 
every step as long as we have some uniform equivalence.

\begin{ass}\label{As-2} The norms $\norms{\cdot}{k}$ and  
$\norm{\cdot}{}$ on $X$ are uniformly equivalent, i.e. there is some 
$C\geq1$ such that 
\begin{equation}
 \frac1C\norm{\cdot}{}\leq \norms{\cdot}{k} \leq C\norm{\cdot}{} 
 \qmz{for all} k\in\N\,.
\end{equation}
In practice $\norms{\cdot}{k}$ is related to the Hessian of some
smooth function at iteration point $x_k$ and $C$ is some 
(usually not explicitly known) bound of that Hessian.
\end{ass}

Notice that the definition of $\gdD{f}{}{x}{h}$ and of 
$\gG{f}{x}{}$ as subset of
$X^*$ merely uses convergence in $X$ and, thus, does not depend on equivalent
norms on $X$. However the Riesz mapping identifying $X^*$ with the Hilbert
space $X$ depends on
the norm. Therefore $\gG{f}{x}{}$ depends on the norm if it is considered as
subset of $X$, which we usually do for simplification of notation. 

\begin{rem}\label{k-norm}
The gradient $\gG{f}{x}{\ep}$ based on the norm  
$\norms{\cdot}{k}$ is understood as subset of $X$ equipped with  
$\norms{\cdot}{k}$ where, in particular, $\cball{x}{\ep}{X}$ is taken 
with respect to $\norms{\cdot}{k}$.
\end{rem}

\subsection{Algorithm}\label{SecAlg1}

\begin{figure}[h!]
\begin{center}
\setlength{\unitlength}{1em}
\setlength{\fboxsep}{0.3em}       
\definecolor{mygray1}{gray}{0.9}
\definecolor{mygray2}{gray}{0.8}
\fns
\begin{picture}(38,35)
\put(0,34.5){\fbox{{\it Initialization:}\: 
    fix parameters $\de$, $\de'$, $\ep_0$,
    functions $T_1$, $T_2$, $G$, initial point $x_0$, and $k:=0$}}

{\thicklines \put(8.5,33.9){\vector(0,-1){1}} }

\put(0,31.8){\fbox{choose a norm $\|\cdot\|_k$, some $a_k\in\pa f(x_k)$ (w.r.t.
    $\|\cdot\|_k$), some $\ep_{k,0}\ge G(\|a_k\|_k,\ep_k)$, and $i:=0$ }}

\put(3.5,5.5){\fcolorbox{black}{mygray2}{
    \makebox(31.2,23.3)[tr]{\bmp{8.5em} {\it Algorithm \ref{A-IA}} \\[-2pt]
     {\tiny (preconditions step size)} \emp}}}

{\thicklines \put(8.5,31.1){\vector(0,-1){2.6}} }

\put(6,27.4){\fcolorbox{black}{white}{$j:=0$,\; choose $a'_0\in\pa f(x_k)$ 
    (w.r.t. $\|\cdot\|_k$, e.g. $a'_0:=a_k$)}}

{\thicklines \put(8.5,26.8){\vector(0,-1){2.5}} }

\put(4.1,23.3){\fcolorbox{black}{white}
    {\bmp{11.6em} \centering {\it null step}\;
         {\tiny $\big($i.e. $0\in\pa^{\ep_{k,i}}\!f(x_k)\big)$} \\[2pt]
           $\|a'_j\|_k\le T_1(\ep_{k,i})$ \emp}}

{\thicklines \put(14.7,22.5){\vector(1,0){2.4}}} \put(15,23){NO}

\put(17.3,23){\fcolorbox{black}{white}
     {\bmp{16em} {\it sufficient descent}\;  
                  $\big(\,h'_j:=\frac{a'_j}{\|a'_j\|_k}\,\big)$ \\[4pt]
      $f(x_k-\ep_{k,i}h'_j)-f(x_k)\le -\de\|a'_j\|_k\ep_{k,i}$ \emp}}

{\thicklines \put(8.5,21.2){\vector(0,-1){3}}} \put(8.7,19.5){YES}
{\thicklines \put(19,20.6){\vector(0,-1){2.5}}} \put(19.2,19.2){YES}
{\thicklines \put(27,20.6){\vector(0,-1){2.6}}} \put(27.2,19.2){NO}

\put(5,17){\fcolorbox{black}{white}
          {\bmp{7em} \centering \reff{L-NSAs} holds with \\[-1pt]
                      $a_{k,i}:=a'_j$ \emp}}

\put(13.4,17){\fcolorbox{black}{white}
          {\bmp{8em} \centering \reff{L-DSAs} holds with \\[-1pt] 
           $a_{k,i}:=a'_j$ $h_{k,i}:=h'_j$ \emp}}

\put(21.7,17){\fcolorbox{black}{mygray1}
         {\bmp{13em} {\it Algorithm \ref{A-IS}} \\
           compute $b'_j\in\pa^{\ep_{k,i}}\!f(x_k)$ with \\[4pt]
           \hs{5} $\dl a'_j,b'_j\dr\le\de'\|a'_j\|_k^2$ \emp}}

{\thicklines \put(27,13.7){\vector(0,-1){1.9}}} 

\put(19.5,10.7){\fcolorbox{black}{white}
  {\bmp{17.5em} with inner approximation of $\pa^{\ep_{k,i}}\!f(x_k)$\\[3pt]
  \hs{1} $A'_j:=\op{conv}
  \big(\{a'_l\,|\, l\le j\}\cup\{b'_l\,|\,l\le j\}\big)$\\[3pt]
  improve norm-minimal element $a'_j$ by \\[3pt]
  \hs{3}  $a'_{j+1}:=\op{arg\; min} \big\{\|a'\|_k\,|\,a'\in A'_j \big\}$ 
\emp}}

{\thicklines \put(34.2,11.8){\line(0,1){13.5}}   
             \put(34.2,25.3){\line(-1,0){22}}     
             \put(12.2,25.3){\vector(0,-1){1} }} 
\put(15,25.7){$j:=j+1$}

{\thicklines \put(2,4.1){\line(0,1){23.5}}     
             \put(2,27.6){\vector(1,0){3.8} }} 
\put(2,28){$i:=i+1$}
\put(0.7,29.6){$k:=k+1$}

{\thicklines \put(8.5,15.3){\vector(0,-1){11.2}}} 
{\thicklines \put(17.7,14.2){\vector(0,-1){12.6}}} 

\put(1.5,3){\fcolorbox{black}{white}
          {{\it null step:}\; 
           $\ep_{k,i+1} := T_2(\|a_{k,i}\|_k)<\ep_{k,i}$}}

{\thicklines \put(0.5,1.6){\vector(0,1){29.5}}} 

\put(0,0.5){\fcolorbox{black}{white}
    {\bmp{30.1em} {\it descent step:}\;  line search in direction $-h_{k,i}$
      for step size $\si_k>\ep_k$ and\\
      \hspace*{6em} $x_{k+1}:=x_k-\si_kh_{k,i}$  \emp}}
\end{picture}
\end{center}

\caption{Flow diagram of Algorithm~\ref{A-Main}.}
\label{MAlg-Flussdiagramm}
\end{figure}
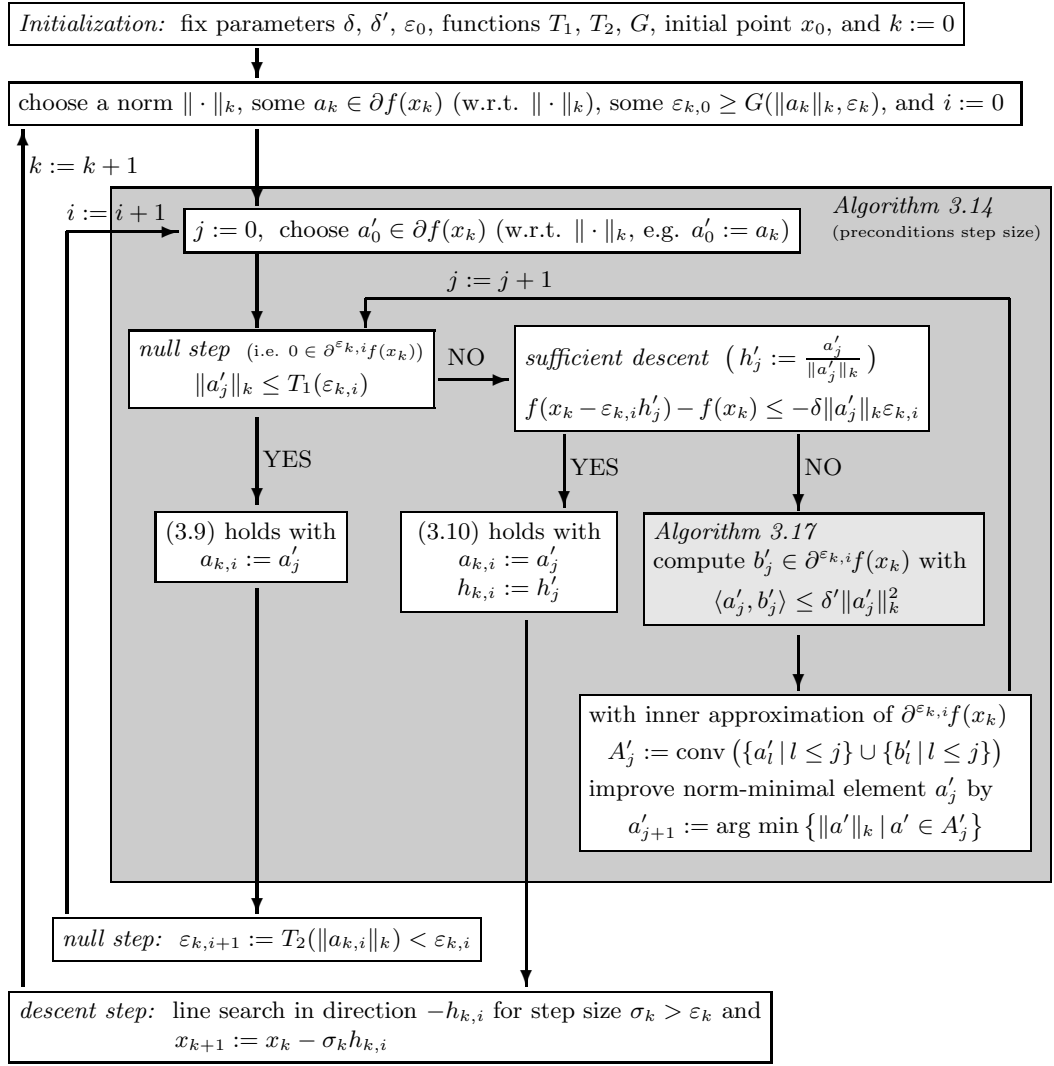

Now we introduce the main algorithm based on two subalgorithms
presented afterwards. We formulate several results that justify the
single steps and that finally show convergence of the algorithm
(cf. Figure \ref{MAlg-Flussdiagramm} below for a rough overview). 
The proofs are collected in Section \ref{Proofs}.

\begin{alg}[Main Algorithm]\label{A-Main}$\;$
\bgl
\item Initialization: Choose $\ThN,\NHI$ and $\USg$ satisfying
   Assumption~\ref{Ass-1},  
\beq 
  \de\in \left]0,1\right[,\;\;\;\;\;\; x_0\in X,\;\;\;\;\;\;
  \ep_0 > 0
\ee
and set $i=k=0$. 
  
\item\label{E-AM-2} Choose some norm $\norms{\cdot}{k}$ subject to  
Assumption \ref{As-2}, some $a_k\in \gG{f}{x_{k}}{}$
(w.r.t. $\norms{\cdot}{k}$), and some
$\ep_{k,0}\geq \USg(\norms{a_k}{k},\ep_k)$.

\item\label{E-AM-3} Determine
$a_{k,i} \in\gG{f}{x_k}{\ep_{k,i}}$ (w.r.t. $\norms{\cdot}{k}$) 
by Algorithm~\ref{A-IA} such that the null step condition
\begin{equation}\label{L-NSAs}
 \norms{a_{k,i}}{k}<\ThN(\ep_{k,i})
\end{equation}
or the condition of sufficient descent   
\begin{equation}\label{L-DSAs}
f\left(x_k-\ep_{k,i}h_{k,i}\right)-f(x_k)\leq
-\de\norms{a_{k,i}}{k}\ep_{k,i}  \qmq{with}
h_{k,i} :=\frac{a_{k,i}}{\norms{a_{k,i}}{k}}
\end{equation}
is satisfied (recall Remark \ref{k-norm} for the meaning of 
$\gG{f}{x_k}{\ep_{k,i}}$ related to $\norms{\cdot}{k}$ and notice that
\reff{L-NSAs} and \reff{L-DSAs} can be satisfied simultaneously).

\item\label{S5SA} In case (\ref{L-NSAs}) set $\ep_{k,i+1}:=\NHI(\ep_{k,i})$,
increment $i$ by one, and go to Step~\ref{E-AM-3}. 

\item \label{descent} If \reff{L-NSAs} is not true, 
choose $\si_k\geq\ep_{k,i}$ such that the condition of sufficient descent
\begin{equation}\label{E-MA-4}
f(x_k-\si_k h_{k,i})-f(x_k)\leq-\delta\norms{a_{k,i}}{k} \si_{k}\,
\end{equation}
is satisfied (notice that $\si_k=\ep_{k,i}$ is always possible,  
since \reff{L-DSAs} is satisfied in this case). 
Then fix the new iteration point 
\begin{equation}\label{iteration}
  x_{k+1}:=x_k-\si_k h_{k,i}\,,
\end{equation}
set $\ep_{k+1}:=\ep_{k,i}$,
increment $k$ by one, set $i=0$, and go to Step~\ref{E-AM-2}.
\el
\end{alg}

\begin{rem}\label{rem1}\leer
\bgl
\item
Instead of $\ep_{k,0}\geq \USg(\norms{a_k}{k},\ep_k)$ in Step~\ref{E-AM-2}
one could also choose 
$$
  \ep_{k,0}\geq \USg(\norms{ a_{k-1}}{k-1},\ep_k) \qmq{for} k>0\,.  
$$

\item The selection of $\si_k$ in Step~\ref{descent} can be done by some 
line search in direction $h_{k,i}$ (cf. Pytlak~\cite{Pyt2009}).

\item One can easily ensure that $\si_k\to 0$ by requiring that e.g. 
\beq
  \ep_k\leq \si_k=\norm{x_k-x_{k+1}}\leq\fkt(\ep_k)
\ee
for some $\fkt:\R_{>0}\to\R_{>0}$ with $\fkt(x)>x$ and
$\lim\limits_{t\to0}\fkt(t)=0,$ since the proof of Theorem~\ref{T-Main} 
shows that $\ep_k\to 0$. But in practice this is usually not necessary.
\el
\end{rem}

The essential point in Algorithm \ref{A-Main} is Step 3 with the computation
of a suitable approximation $a_{k,i}$ of the norm-smallest element 
$\ti a_{k,i}\in\gG{f}{x_k}{\ep_{k,i}}$ (cf. \reff{deff_0}) such that 
the null step condition or the condition of sufficient descent is satisfied
for given $\ep_{k,i}$. Let us briefly discuss the main idea before we
formulate the corresponding subalgorithm. Usually the sets $\gG{f}{y}{}$ 
defining $\gG{f}{x_k}{\ep_{k,i}}$ are not known explicitely. For the algorithm 
we merely suppose that always at least one element $a\in\pa f(y)$ can be
determined numerically 
(cf. Remark \ref{single-el} below for a brief discussion of that point).  
On this basis we select step by step  
elements $b'_j\in\gG{f}{y_j}{}$ for suitable $y_j\in\cball{x_k}{\ep_{k,i}}{X}$
and determine certain $a'_j\in \gG{f}{x_k}{\ep_{k,i}}$ such that, roughly
speaking, the convex hull $A'_j$ of all $a'_l$, $b'_l$ with $l\le j$ is an
approximating subset of $\gG{f}{x_k}{\ep_{k,i}}$ and 
$a'_j$ is a norm-minimal element in $A'_{j-1}$. In doing so we  
still manage that $\|a'_j\|$ is decreasing sufficiently. Therefore 
we reach for $a_{k,i}:=a'_j$ and $j$ large that the null step condition 
\reff{L-NSAs} is satisfied if $0\in\gG{f}{x_k}{\ep_{k,i}}$ or, otherwise,
that $\|a_{k,i}\|$ approximates $\|\ti a_{k,i}\|$ sufficiently well in the
sense of Corollary~\ref{Kiwiel-GS-Verfahren-Corollary}. In the last case
$-h_{k,i}:=-\frac{a_{k,i}}{\norm{a_{k,i}}{}}$ is a descent direction on 
$\cball{x_k}{\ep_{k,i}}{X}$ and, by Proposition \ref{gg-s1} \reff{gg-s1-4},
\begin{equation*}
  f(x_k-\ep_{k,i} h_{k,i}) - f(x_k) \leq 
  \ep_{k,i}\gdD{f}{\ep_{k,i}}{x}{-h} \leq - \de\ep_{k,i}\norm{\ti a_{k,i}} 
\end{equation*}
with $\de\in\,]0,1[$ from Algorithm \ref{A-Main}, i.e. condition \reff{L-DSAs}
of sufficient descent is satisfied with the standard norm. Clearly the quality
of the algorithm is closely related to the quality of the approximating set
$A'_j$ and, in some applications, we can improve the quality substantially 
by choosing a suitable equivalent norm $\norms{\cdot}{k}$ in each step.

Let us now provide the precise algorithm where quantities determined here are
marked by~$'$. 
\begin{alg}\label{A-IA}
Let $T_1$, $\de\in]0,1[$, $x_k\in X$, $\ep_{k,i}>0$, and $\norms{\cdot}{k}$
be as in Step~\ref{E-AM-3} of Algorithm~\ref{A-Main}
for some $k,i\in\N$.
\bgl
\item Choose some $a'_0\in\gG{ f}{x_k}{}$ (w.r.t. $\norms{\cdot}{k}$)
and some $\de'\in\,]\de,1[$ and set $j=0$. 
(Typically, but not necessarily, $a'_0$ agrees with $a_k$ from Algorithm 
\ref{A-Main}.)

\item\label{A-IAS2} Set $a_{k,i}:=a'_j$ and 
$h'_j :=\frac{a'_j}{\|a'_j\|_k}$. If $a_{k,i}$ satisfies the null step
condition \reff{L-NSAs} or condition \reff{L-DSAs} of sufficient descent, 
stop and return $a_{k,i}$.

\item Otherwise compute some $b'_j \in\gG{f}{y_j}{}$ 
(w.r.t. $\norms{\cdot}{k}$) 
for some $y_j \in\big[x_k,x_k-\ep_{k,i}h'_j\big]$ 
by Algorithm~\ref{A-IS} such that 
\begin{equation}\label{a_jb_jkl}
 \sk{ a'_j}{b'_j}{k}\leq \de' \norms{a'_j}{k}^2\,.
\end{equation}

\item Choose some subset 
$B'_j\subseteq\{a'_l\tme  l\leq j\}\cup\set{b'_l}{ l\leq j }$
such that $a'_j,b'_j\in B'_j$ and set 
\begin{equation*}
 A'_j:=\op{conv}\;B'_j\,.
\end{equation*}
\item\label{A-IAS5} Compute 
$$ a'_{j+1}:=\arg\min\set{\norms{a'}{k}}{  a'\in A'_j }\,,
$$
increment $j$ by one, and go to Step \ref{A-IAS2}.
\el
\end{alg}
\noindent
Notice that $a'_l,b'_l\in \gG{f}{x_k}{\ep_{k,i}}$ for all 
$1\le l\le j$ by induction and that
$$ A'_j=\op{conv} B'_j\subseteq\gG{f}{x_k}{\ep_{k,i}}\,.
$$ 
Hence $A'_j$ can be considered as some inner approximation of
$\gG{f}{x_k}{\ep_{k,i}}$ and the norm-smallest element $a'_{j+1}\in A'_j$ 
is an approximation of the norm-smallest element 
$\ti a_{k,i}\in \gG{f}{x_k}{\ep_{k,i}}$. Algorithm~\ref{A-IA} ensures 
with \reff{a_jb_jkl} that $\|a'_j\|$ decreases sufficiently, i.e. 
we have $\|a_{j+1}'\|_k\le\ga'\|a_j'\|_k$ for some $\ga'\in(0,1)$
as long as the null step condition \reff{L-NSAs} is not fulfilled
(cf. the proof of Lemma~\ref{L-IA}). Hence,  
the null step condition {\rm (\ref{L-NSAs})} has to be satisfied for 
some $a_{k,i}=a'_j$ after finitely many steps if we do not meet 
condition \reff{L-DSAs} of sufficient descent before.
In practice we usually take 
\ba
B'_j
&=&
\{a'_0\}\cup\set{a'_l}{ j-m\leq l\leq j}\cup\sset{b'_j} \qmq{or} \nn\\
B'_j
&=&
\{a'_0, a'_j\}\cup\set{b'_l}{ j-m\leq l\leq j }  \nn
\ea
with $m\approx 10$. 

\begin{rem}
Note that the computation of $a'_{j+1}$ is equivalent to the minimization 
of a quadratic function defined on some $\#B'_j$-simplex. This can be
easily done with SQP or semi smooth Newton methods 
(cf. \cite{Alt2002,NocedalWright,Ulbricht2006}). 
Since $\mbox{dim}\,X\gg m$ for typical applications, we can neglect the 
computational time for $a'_{j+1}$ compared to that needed for the computation
of a gradient of $f$. 
\end{rem}

We complete our algorithm with the precise scheme about the selection of 
$y_j$ in Step~3 of Algorithm~\ref{A-IA}
by some nesting procedure for the segment 
$\big[x_k,x_k-\ep_{k,i}h'_j\big]$. New quantities determined in the
subalgorithm are marked by $''$. 

\begin{alg}\label{A-IS} 
Let $0<\de<\de'<1$, $x_k\in X$, $\ep_{k,i}>0$, $\norms{\cdot}{k}$,  
$a'_j\in\gG{f}{x_k}{\ep_{k,i}}$, and $h'_j$  
be as in Step 3 of Algorithm~\ref{A-IA}.
(Notice that both the null step condition \reff{L-NSAs} and condition 
\reff{L-DSAs} of sufficient descent are violated for $a_{k,i}=a'_j$.)   
\bgl
\item\label{E-IS-1} Set $l=0$, $x''_0 :=x_k$, and 
$y''_0 :=x_k-2\ep_{k,i}h'_j$ (notice that
$\frac{x''_0+y''_0}{2}=x_k-\ep_{k,i}h'_j$).
\item\label{E-IS-4} Choose some 
$b''_l \in\gG{f}{\frac{x''_l+y''_l}2}{}$ (w.r.t. $\norms{\cdot}{k}$). 
\item\label{E-IS-2} If $b'_j:=b''_l$ satisfies (\ref{a_jb_jkl}) stop and return
$b'_j$.  
\item\label{E-IS-5} 
Otherwise choose $x''_{l+1}\in\sset{x''_l\,,\,\frac{x''_l+y''_l}2}$ and
$y''_{l+1}\in\sset{\frac{x''_l+y''_l}2\,,\,y''_l}$ such that 
\begin{eqnarray}
 \norms{y''_{l+1}-x''_{l+1}}{k}&=&\tfrac12\norms{y''_{l}-x''_{l}}{k}
\qquad\mbox{and}\nn\\[2mm]  
f(y''_{l+1})-f(x''_{l+1})&>& -\de \norms{a'_j}{k}\norms{y''_{l+1} -
  x''_{l+1}}{k} \,\label{E-IS-6} 
 \end{eqnarray}
where we take $x''_1=x''_0$ if $l=0$
(this way the condition of sufficient descent is
violated on segment $[x''_{l+1},y''_{l+1}]$ with
$a=a'_j$).
\item Increment $l$ by $1$ and go to Step~\ref{E-IS-4}.
\el
\end{alg}
\noindent
A slightly simplified survey about the complete algorithm is given in 
Figure~\ref{MAlg-Flussdiagramm}.

\begin{rem}\label{single-el}
While the implementation of the most steps in Algorithm \ref{A-Main}
and its subalgorithms should be quite
clear, let us briefly discuss how to choose some 
element $a\in\pa f(x)$.
In our applications we usually have a representation of  
$\partial f(x)$ that allows the numerical computation of 
some element $a\in \partial f(x)$. More precisely, in many cases 
$f$ is continuously
differentiable on an open set $U\subset X=\R^n$ such that 
$\R^n\setminus U$ has zero Lebesgue measure. Here we can use 
Proposition~2.1.5 or Theorem~2.5.1 from Clarke \cite{Clarke} 
to get single elements $a\in \partial f(x)$. If $f$
is defined to be the pointwise maximum or minimum of smooth 
functions, Proposition~2.3.12 or Proposition~2.7.3 in \cite{Clarke} can be
used to determine 
some $a\in \partial f(x)$. Moreover we can combine this with other calculus
rules as e.g. the chain rule \cite[Theorem~2.3.9]{Clarke}.
Beyond these methods, that are sufficient for the benchmark problems 
considered in Section~\ref{Benchmark}, also discrete
approximations of elements of $\partial f(x)$ like e.g. in
\cite{Bagirov2008} can be used. Let us finally state that the presented
algorithm assumes the possibility to compute at least one element of 
$\partial f(x)$. 
\end{rem}

Let us now justify the essential steps of the algorithm, i.e. that the 
required conditions can be reached and that the iterations typically  
terminate after finitely many steps. We start with 
Algorithm~\ref{A-IS} and consider in particular Step~\ref{E-IS-5}.

\begin{prop}[properties of Algorithm~\ref{A-IS}]
\label{ConvInterv}
Let the assumptions of Algorithm~{\rm \ref{A-IS}} be satisfied. Then:
\bgl
\item The choice in Step~{\rm \ref{E-IS-5}} of Algorithm~{\rm \ref{A-IS}}
is possible for every $l\in\N$. 
\item \label{ConvInterv-2}
The set 
\begin{equation*}
 \La_l :=\set{\la\in[0,1]}{\mbox{there is some}\;b'_j\in
 \gG{f}{\la x''_l+(1-\la)y''_l}{}\;\mbox{satisfying \reff{a_jb_jkl}\,}} 
\end{equation*}
has positive Lebesgue measure for every $l\in\N$. 
 \item  If Algorithm~{\rm \ref{A-IS}} does not terminate and, therefore,
   produces 
sequences $(x''_l)$ and $(y''_l)$ converging to some
$y''_{\infty}\in[x''_1,y''_1]$, then there is some 
$b'_j\in \gG{f}{y''_{\infty}}{}$ satisfying \reff{a_jb_jkl} and 
$f$ is not strictly differentiable at $y''_{\infty}$. 
\item If $f$ is convex on a neighborhood of $[x''_1,y''_1]$, then
Algorithm~{\rm \ref{A-IS}} terminates in 
Step~{\rm \ref{E-IS-2}} already for $l=0$. 
\el
\end{prop}
Though it is not trusted that we find some $b'_j$ satisfying \reff{a_jb_jkl}
after finitely many steps, there is an extremely good chance according to 
Proposition \ref{ConvInterv} \reff{ConvInterv-2}. In practice 
Algorithm~\ref{A-IS} always terminated, also for rather complex simulations
presented here. Nevertheless there are examples where the algorithm  
does not terminate (at least theoretically) as a simple induction argument 
shows, e.g., for $f(t)=-t^2\sin(2\pi/t)$
with $\ep_{k,i}=1=y''_1=a'_j$ and $f(0)=0=x''_1$.
\begin{rem}
Typically it is much cheaper for time consumption to compute merely
the scalar $\sk{b'_j}{a'_j}{k}$ in (\ref{a_jb_jkl}) instead of the complete
vector $b'_j$. Therefore we compute $b'_j$ only if (\ref{a_jb_jkl}) is
satisfied. 
\end{rem}

\begin{prop}[properties of Algorithm~\ref{A-IA}]\label{SecAlg-s1}
Let the assumptions of Algorithm~{\rm \ref{A-IA}} be satisfied, let $f$ be
Lipschitz continuous on some neighborhood of $\ol{B_{\ep_{k,i}}(x_k)}$, and
suppose that Algorithm~{\rm \ref{A-IS}} always terminates.
Then Algorithm~{\rm \ref{A-IA}} stops after finitely many steps and returns
some $a_{k,i}\in\gG{f}{x_k}{\ep_{k,i}}$ satisfying {\rm (\ref{L-NSAs})} or
{\rm (\ref{L-DSAs})}.
\end{prop}

\begin{prop}[properties of Algorithm~\ref{A-Main}]\label{P-Main}
Let the assumptions of Algorithm~{\rm \ref{A-Main}} be satisfied 
and let $x_k$ be an iteration point from that algorithm. Then: 
\bgl
\item \label{P-Main-1}
If $\ep_{k,i}>0$ is related to $x_k$ for some $i\in\N$,
then there exists $a_{k,i}\in\gG{f}{x_k}{\ep_{k,i}}$ 
satisfying the null step condition {\rm (\ref{L-NSAs})} or condition
{\rm (\ref{L-DSAs})} of sufficient descent.

\item If $0\notin\gG{f}{x_k}{}$, then there are only
  finitely many $i\in\N$ such that {\rm (\ref{L-NSAs})} is satisfied.  
\el
\end{prop}

Though Proposition \ref{P-Main} \reff{P-Main-1} already follows
from Propositions \ref{SecAlg-s1} we will still give a brief 
independent proof of it in the next section. 

Summarizing we can say that, in principle, the presented algorithm always
works and cannot get stuck, i.e. at most finitely many subiterations are 
necessary to find a new iteration point $x_{k+1}$. 
The only point is that Algorithm \ref{A-IS}
might not terminate which, however, is quite unlikely according to 
Proposition \ref{ConvInterv} \reff{ConvInterv-2} and which never happened in
our simulations. 

Let us finally confirm that the presented descent algorithm 
can reach both minimizers and critical points of $f$.

\begin{theo}[accumulation points are critical points]\label{T-Main}
Let the assumptions of Algorithm~{\rm \ref{A-Main}} be satisfied and let 
$(x_k)$ be a corresponding sequence of iteration points. Then
$\big(f(x_k)\big)$ is strictly decreasing. Moreover, if $x$ is an accumulation
point of $(x_k)$, then $0\in\gG{f}{x}{}$ and $f(x_k)\to f(x)$.
\end{theo}

As consequence we can formulate some more precise statement.

\begin{prop}\label{konvergentertheoretischerTheorem2}
Let the assumptions of Algorithm~\ref{A-Main} be satisfied, let 
$(x_k)$ be a corresponding sequence of iteration points with step sizes 
$\si_k\to 0$, and suppose that $\set{x_k}{ k\in\N }$ is relatively compact.
Then each accumulation point of $(x_k)$ is a critical point of $f$
and, if $\set{ x\in X}{ f(x)\leq f(x_0) }$ contains only finitely many
critical points, $(x_k)$ converges to a critical point of $f$.
Moreover, if $(x_k)$ is not convergent, then $(x_k)$ has no isolated
accumulation point. 
\end{prop}

\begin{rem}
If $X=\R^n$ and $\set{ x\in X}{ f(x)\leq f(x_0) }$ is bounded, then
$\set{x_k}{ k\in\N }$ is relatively compact. 
\end{rem}

\subsection{Proofs}\label{Proofs}
\begin{pf}{ of Proposition~\ref{ConvInterv}}
Since $a_{k,i}:=a'_j$ does not satisfy (\ref{L-NSAs}), 
we have $\|a'_j\|_k\neq0$. By construction 
$\|y''_l-x''_l\|_k=\left(\frac1{2}\right)^{l-1}\ep_{k,i}>0$. 

(1) Since (\ref{L-DSAs}) is not fulfilled, we have (\ref{E-IS-6}) for $l=0$ 
with $y_1''=x_k-\ep_{k,i}h'_j$. 
Assume that (\ref{E-IS-6}) holds for $l-1$. Then
\begin{eqnarray}
 f(x''_l)-f\left(\frac{x''_l+y''_l}2\right)+
 f\left(\frac{x''_l+y''_l}2\right)-f(x''_l)>-\de 
 \norms{a'_j}{k}\norms{x''_{l}-y''_{l}}{k}\,.\label{E-IS-8} 
\end{eqnarray}
If neither $x''_{l+1}=x''_l$ and $y''_{l+1}=\frac{x''_l+y''_l}2$ nor 
$x''_{l+1}=\frac{x''_l+y''_l}2$ and $y''_{l+1}=y''_l$ satisfy (\ref{E-IS-6}) 
for $l$, then
\bee
f(x''_l)-f\left(\frac{x''_l+y''_l}2\right)+
f\left(\frac{x''_l+y''_l}2\right)-f(y''_l)\leq 
-2\de\norms{a'_j}{k}\frac12\norms{x''_{l}-y''_{l}}{k}\, 
\ee
which contradicts (\ref{E-IS-8}). Hence the choice of 
$x''_{l+1}$ and $y''_{l+1}$ is always possible and
induction gives (\ref{E-IS-6}) for all $l>0$.

(2) Consider $\tilde{f}_l:[0,1]\to\R$ with
$\tilde{f}_l(t):=f(\,y''_l+t(x''_l-y''_l)\,)$. As Lipschitz continuous
function, $\ti f_l$ is absolutely continuous and differentiable almost
everywhere. Therefore
\bee
\int\limits_0^1\tilde{f}_l'(t)dt= f(x''_l)-f(y''_l)
\stackrel{(\ref{E-IS-6})}{>}-\de \norms{a'_j}{k}\norms{x''_l-y''_l}{k}
\stackrel{\de<\de'}{>} -\de' \norms{a'_j}{k}\norms{x''_l-y''_l}{k}\,.
\ee 
Hence the set of all $t\in[0,1]$ with $\tilde{f}_l'(t)>-\de'
\|a'_j\|_k\norms{x''_l-y''_l}{k}$ has positive measure. 
Using Clarke's chain rule \cite[Theorem~$2.3.10$]{Clarke} and
$x''_l-y''_l=-\|x''_l-y''_l\|_k\, \frac{a'_j}{\|a'_j\|_k}$, we get 
\begin{equation*}
 \tilde{f}'(t)\in
 \big\dl \gG{f}{y''_l+t(x''_l-y''_l)}{},x''_l-y''_l\big\dr_k =
 \frac{-\norms{x''_l-y''_l}{k}}{\|a'_j\|_k}
 \big\dl \gG{f}{y''_l+t(x''_l-y''_l)}{},a'_j\big\dr_k\, 
\end{equation*} 
which implies that $\La_l$ has positive Lebesgue measure.

(3) By Lebourg's mean value theorem (cf.~\cite[Prop.~$2.3.7$]{Clarke}),
there exists some $\tilde{x}''_l\in[x''_l,y''_l]$ and some
$\tilde{b}''_l\in\gG{f}{\tilde{x}''_l}{}$ with 
\begin{equation}\label{E-IS-7}
 -\de\norms{a'_j}{k}\norms{y''_l-x''_l}{k}\stackrel{(\ref{E-IS-6})}{<}
 f(x''_l)-f(y''_l)=\sk{\tilde{b}''_l}{x''_l-y''_l}{k}= 
 -\frac{\norms{y''_l-x''_l}{k}}{\|a'_j\|_k}\sk{\tilde{b}''_l}{a'_j}{k}\,. 
\end{equation}
Clearly $\frac{y''_l+x''_l}2\to y''_{\infty}$ and $\tilde{x}''_l\to
y''_{\infty}$. Since $f$ is Lipschitz continuous near $y''_{\infty}$, 
the sequences $(b''_l)$ and $(\tilde{b}''_l)$ are bounded in the Hilbert space
$X$. Therefore, up to a subsequence, $b''_l\wto:b''$ and 
$\tilde{b}''_l\wto:b'_j$. The upper semicontinuity of Clarke's 
generalized gradient implies $b'',\,b'_j\in\gG{f}{y''_{\infty}}{}$
(cf. \cite[Prop.~$2.1.5$]{Clarke}). 
Since we assume that (\ref{a_jb_jkl}) always fails in Step 3 of 
Algorithm~\ref{A-IS}, we have  
\begin{equation}
 \sk{\tilde{b}''_l}{a'_j}{k}\stackrel{(\ref{E-IS-7})}{<}\de\norms{a'_j}{k}^2\;
 \stackrel{\de<\de'}{<}\; \de'\norms{a'_j}{k}^2 \;<\;\sk{b''_l}{a'_j}{k}\;. 
\end{equation}
Taking the limit we obtain that $b'_j$ satisfies (\ref{a_jb_jkl}) 
and that $b'_j\neq b''$. Hence $\gG{f}{y''_{\infty}}{}$ is not a singleton
and, consequently, $f$ is not strictly differentiable at $y''_{\infty}$ 
(cf. \cite[Prop.~$2.2.2$]{Clarke}).

(4) Recall that Clarke's generalized gradient $\gG{f}{x}{}$ agrees with the 
(convex) subdifferential for a convex $f$ (cf. \cite[Prop.~$2.2.7$]{Clarke}). 
Since (\ref{L-DSAs}) is not satisfied for $a_{k,i}=a'_j\in\gG{f}{x''_0}{}$,
the definition of the subdifferential gives for any
$b''_0\in \pa f\big(\frac{x''_0+y''_0}{2}\big)$ that
\begin{equation*}
 \de' \ep_{k,i}\norms{a'_j}{k}       \stackrel{\de<\de'}{>}     
 \de \ep_{k,i}\norms{a'_j}{k}     >    
 f(x''_0) - f\Big(\frac{x''_0+y''_0}{2}\Big)    \geq 
 \big\dl b''_0 , \ep_{k,i}\tfrac{a'_j}{\|a'_j\|_k}\big\dr_k 
\end{equation*}
which implies (\ref{a_jb_jkl}) for $b'_j=b''_0$.
\end{pf}

\medskip

As preparation for the proof of Proposition \ref{SecAlg-s1}
we consider some technical lemma. Notice that the statement remains true with
$\|\cdot\|_k$ instead of $\|\cdot\|$.

\begin{lem}\label{L-IA}
Let $\ga\in\left[0,1\right[$ and $L>0$ be constants and
let $(a'_j)$ and $(b'_j)$ be sequences in $X$ 
such that for all $j\in\N$: 
\bgl
\item\label{L-IAE1} $\sk{ a'_j}{ b'_j}{}\leq \ga\|a'_j\|^2$,
\item\label{L-IAE2} $\|b'_j\|\leq L $\,,
\item $\|a'_{j+1}\|\leq
\min\limits_{\la\in\left[0,1\right]}\norm{\la a'_j+(1-\la)b'_j }$\,.
\el
Then $\norm{ a'_j}\to 0$.
\end{lem}

\begin{pf}{}
The sequence $\big(\norm{a'_j}{}\big)$ 
is non increasing by (3) and, hence, convergent.
Since $L$ is just a bound for the $b'_j$ according to (2), we can assume that 
$\norm{a'_j}{}<L$ for all $j\in\N$. Therefore
$(1-\ga)\norm{a'_j}{}^2<L^2$ and
\begin{equation}\label{laj}
 \la_j:=1-\frac {(1-\ga)\norm{a'_j}{}^2}{(1-\ga)^2\norm{a'_j}{}^2+L^2} 
 \in\: ]0,1] \,.
\end{equation}
We derive
\begin{eqnarray*}
\norm{ a'_{j+1}}^2&\leq&\norm{\la_ja'_j+ (1-\la_j)b'_j}{}^2\\
         &=
         &\la_j^2\norm{a'_j}{}^2+2\la_j(1-\la_j)\sk{a'_j}{b'_j}{}+(1-\la_j)^2
          \norm{b'_j}{}^2      \\ 
         &\stackrel{\reff{L-IAE1}, \reff{L-IAE2}}{\leq}& 
\big( \la_j^2 +2\la_j(1-\la_j)\ga\big)\norm{a'_j}{}^2 +(1-\la_j)^2L^2\\
         &\leq& \big( \la_j +(1-\la_j)\ga\big)^2\norm{a'_j}{}^2 +(1-\la_j)^2L^2\\
&\stackrel{(\ref{laj})}{=}&
\bigg(\frac{(1-\ga)^2\norm{a'_j}{}^2+L^2-(1-\ga)\norm{a'_j}{}^2+
\ga(1-\ga)\norm{a'_j}{}^2}{(1-\ga)^2\norm{a'_j}{}^2+L^2}\bigg)^2\norm{a'_j}{}^2\\
&&
\; + \;
\frac {(1-\ga)^2\norm{a'_j}{}^4L^2}{\big((1-\ga)^2\norm{a'_j}{}^2+L^2\big)^2}\\
&=&\frac{L^4\norm{a'_j}{}^2}{\big((1-\ga)^2\norm{a'_j}{}^2+L^2\big)^2}+
\frac{L^2(1-\ga)^2\norm{a'_j}{}^4}{\big((1-\ga)^2\norm{a'_j}{}^2+L^2\big)^2}\\
&=&\frac{L^2\norm{a'_j}{}^2}{(1-\ga)^2\norm{a'_j}{}^2+L^2}\,.
\end{eqnarray*}
Taking the limit we obtain with $\al:=\lim_{j\to\infty}\norm{a'_j}{}^2$ that
$\al \leq \frac{\al L^2}{(1-\ga)^2\al+L^2}$.
But this is only possible for $\al=0$, which implies the assertion. 
\end{pf}

\medskip

\begin{pf}{ of Proposition \ref{SecAlg-s1}}
Let $L$ be the Lipschitz constant of $f$ on some neighborhood of
$\cball{x_k}{\ep_{k,i}}{k}$.
Thus the $b'_j$ are bounded by $L$ according to Proposition~\ref{gg-s1}.
With (\ref{a_jb_jkl}) and the definition of $a'_{j+1}$ we verify the
assumptions of Lemma~\ref{L-IA} with $\ga=\de'$ and the norm $\|\cdot\|_k$.
We obtain the claim,
since the algorithm stops as soon as $\|a'_j\|_k\leq \ThN(\ep_{k,i})$. 
\end{pf}

\medskip

\begin{pf}{ of Propostion~\ref{P-Main}}
For the first assertion we fix $k,i\in\N$. If
$0\in\gG{f}{x_k}{\ep_{k,i}}$ (w.r.t. $\norms{\cdot}{k}$), then $a_{k,i}:=0$
satisfies (\ref{L-NSAs}). Otherwise we choose $a_{k,i}:=\ti a\in\pa f(x_k)$
and have $h_{k,i}=-\ti h$ with $\ti a$, $\ti h$ as in
Proposition~\ref{Satz1Umg} for $x=x_k$ and $\ep=\ep_{k,i}$.
By Lebourg's mean value theorem 
(cf.~\cite[Prop.~$2.3.7$]{Clarke}) there are
$y\in[x_k,\,x_k-\ep_{k,i}h_{k,i}]$ and
$a\in\gG{f}{y}{}\subseteq\gG{f}{x_k}{\ep_{k,i}}$ 
(w.r.t. $\norms{\cdot}{k}$) such that 
\begin{eqnarray*}
f(x_k-\ep_{k,i} h_{k,i})-f(x_k)=   \ep_{k,i}\sk{a}{-h_{k,i}}{k} =
\ep_{k,i}\sk{a}{\ti h}{k}
\stackrel{(\ref{supportgG})}{\leq} \ep_{k,i}\,\gdD{f}{\ep_{k,i}}{x}{\ti
  h}\stackrel{(\ref{WeftverR})}{<} -\de\ep_{k,i}\norms{a_{k,i}}{k}\,. 
\end{eqnarray*}
Thus (\ref{L-DSAs}) is satisfied and the claim is verified.
(Notice that the statement also follows from Proposition~\ref{ConvInterv}
and Proposition~\ref{SecAlg-s1} if the corresponding assumptions are
satisfied.)   

For the second assertion we fix $k$ and $\norms{\cdot}{k}$. Since
$0\notin \gG{f}{x_k}{}$, Lemma~\ref{L-BW} provides $\ep$, $K>0$ such that 
\bee
  \norms{a}{k}>K \qmz{for all} a\in \gG{f}{x_k}{\ep}\,.
\ee
If the algorithm would make infinitely many null steps, 
i.e. (\ref{L-NSAs}) holds
for infinitely many $i\in\N$, then $\ep_{k,i}\to 0$ and $\ThN(\ep_{k,i})\to 0$
by $\ep_{k,i+1}=T_2(\ep_{k,i})$ and Assumption~\ref{Ass-1}. Hence 
$\ep_{k,i}<\ep$ for large $i$ and, thus,
\begin{equation*}
 0 < K \leq \norms{a_{k,i}}{k}\stackrel{(\ref{L-NSAs})}{\leq}
 \ThN(\ep_{k,i}) \to 0 \,. 
\end{equation*}
But this is a contradiction and verifies the assertion. 
\end{pf}

\medskip

\begin{pf}{ of Theorem~\ref{T-Main}}
$\big(f(x_k)\big)$ is strictly decreasing by construction
(cf. \reff{E-MA-4}). Let us suppose that $x$ is an accumulation point of 
$(x_k)$. Then $(f(x_k))$ has a unique accumulation
point and, by continuity of $f$, we have $f(x_k)\to f(x)$.

By $i_k$ we denote the index $i$ related to $k$ 
leading to the new iteration point $x_{k+1}$ in Step~\ref{descent} of 
Algorithm~\ref{A-Main}
(notice that the assumptions of the theorem imply the existence of
$i_k$). Then we have $\ep_{k+1}=\ep_{k,i_k}$ and we set $\hat a_k:=a_{k,i_k}$.
Since the null step condition (\ref{L-NSAs}) is never satisfied for
$a_{k,i_k}$ and since $\si_k\ge\ep_{k,i_k}$, we have for all $N\in\N$
\begin{eqnarray*}
f(x)-f(x_0)&\leq&f(x_{N+1})-f(x_0)= \sum\limits_{k=0}^Nf(x_{k+1})-f(x_k)\\
           &\stackrel{(\ref{E-MA-4})}{\leq}& \sum\limits_{k=0}^N-\de
           \si_k\norms{\hat a_k}{k} \leq -\de\sum
           \limits_{k=0}^N\ep_{k+1}\ThN(\ep_{k+1})\,. 
\end{eqnarray*}
Hence the right hand side is bounded independent of $N$. Therefore 
$\ThN(\ep_k)\ep_k\to 0$ and, since $\ThN$ is nondecreasing,   
\begin{equation}\label{epkTONull}
\ep_k\to 0\,.
\end{equation}

For contradiction we assume that $0\not\in \pa f(x)$. 
By Lemma~\ref{L-BW} we find 
$\ep>0$ and $h\in X$ with $\norm{h}=1$ such that 
\bee
  \gdD{f}{\ep}{x}{h}<0 \quad \text{(w.r.t. $\norm{\cdot}$)}.
\ee 
With $C\ge 1$ from Assumption~\ref{As-2} we have for all 
$a\in\gG{f}{x_k}{\ep_{k,i}}$ (w.r.t. $\norms{\cdot}{k}$)
\ben \label{T-Main-p1}
  \norms{a}{k}\stackrel{\norm{h}=1}{\geq} 
  \norms{a}{k}\frac{\norms{-h}{k}}{C}
  \geq -\frac1C\sk{a}{h}{k}
  \stackrel{(\ref{supportgG})}{\geq} 
  - \frac1C \gdD{f}{\ep_{k,i}}{x_k}{h} 
  \osm{\reff{gg-e2}}{=} -\frac1C \;
  \sup\limits_{\substack{y\in X\\\norms{y-x_k}{k}\le \ep_{k,i}}} 
  \gdD{f}{}{y}{h} \,.
\ee
For $k,i\in\N$ with 
\ben \label{T-Main-p2}
  \norm{x_k-x}<\frac{\ep}{2} \qmq{and} 
  \ep_{k,i}<\frac{\ep}{2C}\,
\ee
and all $y\in X$ with $\|y-x_k\|_k\le \ep_{k,i}$ we have 
\bee
  \|y-x\|\le \|y-x_k\| + \|x_k-x\| \le C\|y-x_k\|_k + \frac{\ep}{2} 
  \le C\ep_{k,i} + \frac{\ep}{2} \le \ep\,.
\ee
Therefore 
\bee
  \sup\limits_{y\in X:\; \norms{y-x_k}{k}\le \ep_{k,i}} 
  \gdD{f}{}{y}{h} \le
  \sup\limits_{y\in X:\;\norms{y-x}{}\le\ep} \gdD{f}{}{y}{h} 
  \osm{\reff{gg-e2}}{=} \gdD{f}{\ep}{x}{h}\,.
\ee
Consequently \reff{T-Main-p1} gives
\ben \label{Unterschranke_d}
  \|a\|_k \ge -\frac{1}{C}\gdD{f}{\ep}{x}{h} =: K>0\,
\ee
for all $a\in\gG{f}{x_k}{\ep_{k,i}}$ (w.r.t. $\norms{\cdot}{k}$)
and all $k,i$ satisfying \reff{T-Main-p2}.

\medskip

(i) As a first case we assume that there is some $k_0\in\N$ with
$\norm{x_k-x}<\frac{\ep}{2}$ for all $k\ge k_0$. Recall that
\bee
   \ep_{k,0}\ge G(\|a_k\|_k,\ep_k)\,.
\ee
If $k\ge k_0$ and $i_k>0$, then (\ref{L-NSAs}) is valid for $i=i_k-1$ and,
with \reff{Unterschranke_d},
\ben \label{T-Main-p3}
  T_1(\ep_{k,i_k-1})>K \qmq{if} \ep_{k,i_k-1}<\frac{\ep}{2C}\,.
\ee
By Assumption~\ref{Ass-1} there is $\ep_K>0$ such that $T_1(t)\le K$ for all
$t\le\ep_K$. Therefore 
\bee
  \ep_{k,i_k-1}\ge\min\big\{\ep_K,\tfrac{\ep}{2C}\big\}=:\ti\ep>0
  \qmq{for} k\ge k_0\,
\ee
(otherwise $K\ge T_1(\ep_K)\ge T_1(\ep_{k,i_k-1})>K$ by \reff{T-Main-p3}).
Since $T_2$ is not decreasing, 
\bee
  \ep_{k,i_k}=T_2(\ep_{k,i_k-1})\ge T_2(\ti\ep) \qmq{if} i_k>0 \zmz{and}
  k\ge k_0\,. 
\ee
Thus, for general $i_k$ and $a_k\in\pa f(x_k)$ from Step 2 of 
Algorithm~\ref{A-Main},
\bee
  \ep_{k+1}=\ep_{k,i_k} \ge \min\{\ep_{k,0},T_2(\ti\ep)\} \ge
  \min\{ G(\|a_k\|_k,\ep_k),T_2(\ti\ep)\} \qmq{for} k\ge k_0\,.
\ee
From $\ep_k\to 0$ we obtain $G(\|a_k\|_k,\ep_k)\to 0$. Then we can discuss
the two cases of Assumption~\ref{Ass-1} for $G$ separately:

\bgla
\item
We get $\|a_k\|_k\to 0$. This contradicts
\reff{Unterschranke_d}, since $\ep_{k+1}=\ep_{k,i_k}<\frac{\ep}{2C}$ and
$a_k\in\gG{f}{x_k}{\ep_{k,i_k}}$ for all $k$ large enough.

\item
Using \reff{Unterschranke_d} we have 
$\ep_{k+1}\geq\USg(\|a_k\|_k,\ep_k)\geq \ep_k$ for all $k$ sufficiently
large. But this contradicts $\ep_k\to 0$. 
\el
The contradictions imply that $0\in\pa f(x)$ in the special case (i).
Notice that we can argue analogously in the case of Remark \ref{rem1}.1.

\medskip

(ii) Now we assume that $\norm{x_k-x}\geq \frac{\ep}{2}$ 
for infinitely many $k$.
Since $x$ is an accumulation point of $(x_k)$, we can choose a
subsequence $(x_{k(i)})_i$ such that 
\beq
  \ep_k<\frac{\ep}{2C} \qmz{for all} k>k(0)\,,
\ee
\beq
\norm{x_{k(2j)}-x}<\frac{\ep}{4}\,,\quad 
\norm{x_{k(2j+1)}-x}\geq \frac{\ep}{2}
\ee
and
\beq 
\norm{x_{l}-x}<\frac{\ep}{2}
\qmq{for} k(2j)\leq l<k(2j+1) \,.
\ee
Using $I_N:=\set{ j\in\N}{k(2j+1)<N}$ and $f(x_{k+1})\leq f(x_k)$ we get 
\begin{eqnarray*}
 f(x)-f(x_0)&=&\lim\limits_{N\to\infty}f(x_{N})-f(x_0)\\
         &\leq& \lim\limits_{N\to\infty}\sum\limits_{j\in
           I_N}f(x_{k(2j+1)})-f(x_{k(2j)})\\ 
         &=& \lim\limits_{N\to\infty}\sum\limits_{j\in
           I_N}\sum\limits_{l=k(2j)}^{k(2j+1)-1}f(x_{l+1})-f(x_{l})\\ 
         &\stackrel{(\ref{E-MA-4})}{\leq}&
         \lim\limits_{N\to\infty}\sum\limits_{j\in
           I_N}\sum\limits_{l=k(2j)}^{k(2j+1)-1}-\de\norms{\hat a_l}{l}
           \norms{x_{l+1}-x_l}{l}\\ 
         &\stackrel{(\ref{Unterschranke_d})}{\leq}& -\de
         K\lim\limits_{N\to\infty}\sum\limits_{j\in
           I_N}\sum\limits_{l=k(2j)}^{k(2j+1)-1}\norms{ x_{l+1}-x_l}{l}\\ 
 &\stackrel{\tx{Ass.~}\ref{As-2}}{\leq}& -\de \frac
 KC\lim\limits_{N\to\infty}\sum\limits_{j\in
   I_N}\sum\limits_{l=k(2j)}^{k(2j+1)-1}\norm{ x_{l+1}-x_l}\\ 
         &\leq& -\de \frac KC\lim\limits_{N\to\infty}\sum\limits_{j\in
           I_N}\norm{x_{k(2j+1)}-x_{k(2j)}}\\ 
         &\leq& -\de \frac KC\lim\limits_{N\to\infty}\sum\limits_{j\in
           I_N}\frac{\ep}{4}=-\infty\,
\end{eqnarray*}
which is impossible.

Hence we have a contradiction in both cases (i) and (ii). Therefore 
$0\in\gG{f}{x}{}$ and the assertion is verified. 
\end{pf}

\medskip

\begin{pf}{\;of\;Prop.~\ref{konvergentertheoretischerTheorem2}}\\
By Theorem~\ref{T-Main} accumulation points are critical points. Finite sets
consist only of isolated points. The rest follows from general properties of
accumulation points stated in the subsequent Proposition~\ref{PropAP}. 
\end{pf}

\begin{prop}\label{PropAP}$\;$\\
Let $(x_k)$ be a sequence in the Hilbert space $X$ such that 
$\set{x_k}{ k\in\N }$ is relatively compact
and let $A$ be the set of all accumulation points of $(x_k)$.
Then $A\neq\emptyset$ and
$$
\op{dist}(x_k,A):=
\inf\set{ \norm{y-x_k} }{ y\in A }\to0 \qmq{as} k\to\infty\,.
$$
If in addition $\norm{x_k-x_{k+1}}\to 0$ as $k\to\infty$, then either 
\beq
  A=\sset{x} \qmz{for some} x\in X \quad
  (\mbox{i.e. $x_k\osm{k\to\infty}{\lto}x$})
\ee
or $A$ has no isolated points.
\end{prop}
\begin{pf}{}
Clearly $A\ne\emptyset$ by relative compactness.
Let us now assume that $\op{dist}(x_k,A)\nrightarrow0$. Then there is some $K>0$
and a subsequence $(x_{k(i)})_i$ with $\op{dist}(x_{k(i)},A)>K$ for all
$i\in\N$ which contradicts relative compactness.

Under the additional condition $\norm{x_k-x_{k+1}}\to 0$ we assume that
$A$ has an isolated point $x$, i.e. there is some $r>0$ with
$\ball{x}{r}{}\cap A=\sset{x}$. Since $x$ is accumulation point of $(x_k)$,
there exists some $k_0\in\N$ with 
\beq
  x_{k_0}\in \ball{x}{\frac r4}{} \qmq{and}  
  \norm{x_k-x_{k+1}}<\frac r4  \qmz{for all} k\geq k_0\,.
\ee
Moreover, for a possibly larger $k_0$, we have   
\beq
  \op{dist}(x_k,A)<\frac r4 \qmz{for all} k\geq k_0\,,
\ee
since otherwise there is a subsequence $(x_{k(i)})$ with  
$\op{dist}(x_{k(i)},A)\ge\frac r4$ for all $i$,
which is impossible by relative compactness. We thus conclude that
$x_{k_0+1}\in\ball{x}{\frac r4}{}$ and, by induction,
$x_k\in \ball{x}{\frac r4}{}$ for all $k\geq k_0$. 
Consequently $x$ is the only accumulation point of $(x_k)$. 
\end{pf}

\section{Benchmark problems}\label{Benchmark}

We now consider several classical benchmark problems on $\R^n$ 
and compare our numerical results based on Algorithm~\ref{A-Main}
with that of other algorithms found in the literature. In our algorithm 
we normally choose
\begin{equation}\label{standart-de}
 \de'=0.35,\;\de=0.3\;\;\;\mbox{and}\;\;\;\NHI(x)=0.35\cdot x\,.
\end{equation}
Recall that the norm can be changed in each iteration step. We mainly consider
the following two specializations of Algorithm~\ref{A-Main}
(briefly called Algorithm~\ref{A-Main}.A and \ref{A-Main}.B):
\begin{enumerate}[(A)]
\item\label{Benchbasic} We do not change the norm, i.e. 
$\norms{\cdot}{k}:=\norm{\cdot}$ in every iteration step $k$ 
(typically $\|\cdot\|$ is the Euclidean norm). Furthermore we set
(if nothing else is stated) 
\begin{equation}\label{standart-fkt}
 \USg(x,y)=y\,, \quad \ThN(x)=\frac x{\ep_0}
\end{equation}
where $\ep_0$ is the initial radius of the initial neighborhood
$\ball{x_0}{\ep_0}{}$. 

\item\label{BenchNewton} We adapt the norm
$\norms{\cdot}{k}:=\norms{A_k\cdot}{}$ in each iteration $k$ 
with some suitable symmetric and positively definite matrix $A_k$. Moreover we
choose
\beq
\USg(x,y)=x\,, \quad \ThN(x)<x\,.
\ee 
In Step 2 of Algorithm \ref{A-Main} we first determine some 
$\tilde{a}_k\in\gG{f}{x_k}{}$ with respect to the standard
norm~$\norm{\cdot}$. Then we set   
\beq
  a_k:=A_k^{-2}\tilde{a}_k \qmq{and} \ep_{k,0}:=\norms{a_k}{k}
\ee
where $a_k\in\gG{ f}{x_k}{}$ w.r.t. $\norms{\cdot}{k}$, 
since for every $h\in X$
\beq
\sk{a_k}{h}{k}=\sk{A_ka_k}{A_kh}{}=\sk{A_kA_k^{-2}\tilde{a}_k}{A_kh}{}=
\sk{\tilde{a}_k}{h}{}\,. 
\ee
In Step 5 we take $\si_k\ge\ep_k$ ($\ge\ep_{k,i}$) and 
Algorithm~\ref{A-IA} always starts with $a_0'=a_k$. 
\end{enumerate}

The benchmark functions $f:\R^n\to\R$ considered for minimization 
are of course locally Lipschitz continuous. Some of them are even smooth,
but show similar numerical difficulties as nonsmooth problems due to large
second derivatives. In all cases we have an explicit formula for the
computation of at least one element of the generalized gradient (cf. Remark
\ref{single-el}).

We will compare our numerical results  
with results of the bundle method {\bf (BM)} 
and the bundle trust region method {\bf (BTR)} (cf. Alt \cite{Alt2002,Alt2004}),
since these methods for convex functions inspired to some extend 
the development of our algorithm for (not necessarily convex) locally
Lipschitz continuous functions.
Furthermore we consider the {\bf BFGS} algorithm\footnote{There are
different implementations of BFGS (in Matlab there is e.g. 
{\bf fminunc} and {\bf E04KAF}, cf. \cite{Alt2004}).} 
which has been recently studied e.g. by Lewis and Overton. 
They minimized also nonsmooth functions and got promising results
(cf. \cite{LO08,LO13}). In addition we compare with the gradient sampling
algorithm ({\bf GS}). Since a benchmark function $f$ is typically treated by
different algorithms in the literature, we compare all algorithms for 
$f$ fixed.

\subsection{Wolfe function}
The Wolfe function $f:\R^2\to\R$ is a classical convex
benchmark function given by 
\begin{equation}\label{Wolfe}
f(x,y):=\left\{ \begin{array}{lcl}9x+16|y|-x^9& \qmq{if}&x\leq0\,,\\ 9x+16|y|
    & \qmq{if} &0<x<|y|\,,\\5\sqrt{9x^2+16y^2} & \qmq{if} & |y|\leq
    x\,. \end{array}\right.   \tag{Wolfe} 
\end{equation}
This was one of the first functions showing that, even  
in the simple case of a convex function on the 2-dimensional
Euclidean Hilbert space, the classical steepest
descent algorithm might converge to a point that is different from the unique
minimizer and that is not even a critical point.
One easily shows that $(-1,0)$ is the unique minimizer with value
$f(-1,0)=-8$. 
But steepest descent algorithms often converge to the point $(0,0)$
(cf. \cite{Alt2002,Alt2004}) which cannot be a minimizer or 
critical point by convexity and $f(-1,0)<f(0,0)$. 

We compare various algorithms applied to the Wolfe function in 
Alt \cite{Alt2002,Alt2004} with our results. The starting point is always 
$(5,4)$. We apply Algorithm~\ref{A-Main}.\ref{Benchbasic}
with $\ep_0:=0.9$ and we stop as soon as the 
deviation $f(x_k)-f(-1,0)$ is smaller than $10^{-8}$.  

\bigskip

\begin{tabular}{||l|c|c|c|c||}
\hline\hline
Algorithm     & BM &  BTR &  BFGS & Algorithm~\ref{A-Main}.\ref{Benchbasic}\\ 
\hline\hline
Iterations     &                            & $26$       & $21$       & $16$\\
Gradients      &$37$                        & $37$       &            &$28$\\
$f(x_k)-f(-1,0)$& $\approx 1.4\cdot10^{-10}$ &$ <10^{-9}$ & $<10^{-8}$ &
$\approx 2.9\cdot10^{-12} $\\ 
\hline\hline
\end{tabular}

\bigskip

One can get even better results with other choices of parameters
in \reff{standart-de}, but for comparability we wanted to keep these 
parameters fixed for all the benchmark problems.
Nevertheless Algorithm~\ref{A-Main}.\ref{Benchbasic}
always approximates the minimizer well after
only a few iterations and gradients also with other parameters
(usually 25-35 gradients are needed to get $f(x_k)-f(-1,0)<10^{-8}$). 
Thus it seems that we need less or not more iterations and gradients 
than the bundle methods and comparably many as BFGS. 

\subsection{q-max}
Next we consider the q-max function $f:\R^n\to\R$ given by
\begin{equation}\label{qmax}
f(x):=\max\set{x_i^2}{1\leq i\leq n} \qmz{for}
x=(x_1,\ldots,x_n)\in\R^n\,.\tag{q-max} 
\end{equation}
It was used in Alt \cite{Alt2004} to demonstrate that the bundle method and 
the bundle trust region method are fast and stable algorithms for nonsmooth
convex functions. 
Both algorithms have been applied to the three starting points 
\begin{eqnarray*}
u_+    &:=& (1,2,3,\ldots,n)\;\\
v      &:=&0.1\cdot u_+\;,\\
u_{\pm}&:=& (u_{\pm,1},u_{\pm,2},\ldots,u_{\pm,n})
\end{eqnarray*}
where $u_{\pm,i}:=i$ for $i \leq \frac n2$ and $u_{\pm,i}:=-i$ otherwise.
(In \cite{Alt2004} also starting point $e:=(1,\ldots,1)$ was studied. 
Since $\gG{f}{e}{}$ is not single-valued and $\la e\in \partial f(e)$ for some
$\la\in\R$, an exact linesearch in direction $-\frac{e}{\norm{e}}$
would directly find the  global minimizer and, therefore, making the
minimization trivial. 
Unfortunately we do not know which gradients were chosen in \cite{Alt2004}.)

We apply Algorithm~\ref{A-Main}.\ref{Benchbasic} with $\ep_0=0.5$ and 
the special choice $\ThN(x)=15\frac x{\ep_0}$. We stop the line search at the
first point where the function is not decreasing.\footnote{We approximate
  this point numerically and do not compute it analytically.}
The results for the bundle and the bundle trust region method presented for
comparison are taken from \cite{Alt2004}.

A simple argument shows that we get essentially the same iteration points
for all starting points $u_+,\,u_{\pm}$ and $v$, except for scaling with
$0.1$ in the case of $v$ and changing sign for the second half of the vectors
in the case of $u_\pm$. Therefore the values of $f$ at the iteration points are 
the same, except for multiplying by $0.01$ in the case of $v$. We refrain from
giving a proof for that, since we observed this behavior in our numerical
computations too. 

\bigskip

{\small\begin{tabular}{||l|c c c cc c||}
\hline\hline
&\multicolumn{6}{c||}{$n=20$}\\
Algorithm & \multicolumn{4}{c}{BM} &  \multicolumn{2}{c||}
{Algorithm~\ref{A-Main}.\ref{Benchbasic}}\\ 
\hline
Initial point &\multicolumn{2}{c} { $u_+$ }& \multicolumn{2}{c} { $u_{\pm}$
}&  $u_+,\,u_{\pm}$ & $v$ \\ 
Iterations & \multicolumn{2}{c} {}& \multicolumn{2}{c}
{}&\multicolumn{2}{c||}{$142$}\\ 
Gradients  &\multicolumn{2}{c} {$247$}              &\multicolumn{2}{c}{$199$
}&\multicolumn{2}{c||}{$246$}\\ 
Value     &\multicolumn{2}{c} {$1.090\cdot 10^{-9}$}&\multicolumn{2}{c}
{$4.145\cdot 10^{-9}$}& $1.4\cdot10^{-10}$&$1.4\cdot10^{-12}$\\ 

\hline\hline

&\multicolumn{6}{|c||}{$n=50$}\\
Algorithm & \multicolumn{2}{c}{BM} &  \multicolumn{2}{c}{BTR} &
\multicolumn{2}{c||}{Algorithm~\ref{A-Main}.A}\\
\hline
Initial point & $u_+$ & $v$ &  $u_+$ & $v$ &  $u_+,\,u_{\pm}$ & $v$ \\
Iterations &       &         &                 &                   &\multicolumn{2}{c||}{$126$}\\
Gradients  & $3,108$&$3,140$   &$451$            &$321$              &\multicolumn{2}{c||}{$311$}\\
Value      &$95.11$&$0.01316$&$1.3\cdot 10^{-7}$&$9.6\cdot 10^{-7}$&$9.6\cdot 10^{-6}$&$9.6\cdot 10^{-8}$\\
\hline\hline
\end{tabular}}

\bigskip

\noindent
More iterations of Algorithm~\ref{A-Main}.\ref{Benchbasic} for $n=50$ and 
initial points $u_+$, $u_\pm$ give
\bgl[--]
\item value $1.9\cdot 10^{-9}$ after $175$ iterations and $452$ gradients and
\item value $2.0\cdot 10^{-11}$ after $200$ iterations and $537$ gradients
\el
(for initial point $v$ we get the values $1.9\cdot 10^{-11}$ and $2.0\cdot
10^{-13}$, respectively).
Unfortunately we did not find results with other starting points for further 
comparisons. 

The Newton method is just the steepest descent method with proper step
size. Thus the step size strategy is crucial, as we can also see from the fact
that the bundle trust region method gives much better results than the bundle
method. 
In particular one could expect that Algorithm~\ref{A-Main}.\ref{BenchNewton} 
terminates at the minimizer after at most $n$ steps for every starting point
and properly chosen gradients. In practice we applied
Algorithm~\ref{A-Main}.\ref{BenchNewton} with $\si_k:=\ep_k$ and 
$A_k:=\op{id}$ as approximation of the Hessian (i.e. we always took 
the Euclidian norm $\norms{\cdot}{k}=\norm{\cdot}$) and it really stopped 
exactly at the minimizer after $n$ steps and $n$ gradient
computations for all $3$ initial points and 
both cases $n=20$ and $n=50$.
Notice that the Hessian of $f$ is typically not regular or it is even not
defined. Since we require $A_k$ to be regular, we set always $A_k=\op{id}$
which is in fact a (scaled) mean value of all possible Hessians. 
The crucial difference to Algorithm~\ref{A-Main}.\ref{Benchbasic} is the 
choice of the precise Newton step size $\ep_{k,0}:=\norm{a_k}$.

Summarizing we can say that Algorithm~\ref{A-Main}.\ref{Benchbasic} 
gives a good approximation of the minimizer
after relatively few iterations and gradient computations. Again it appears
that both versions of Algorithm~\ref{A-Main} are faster than 
the bundle methods. Due to the essentially quadratic structure of the
function it is not surprising that our basic algorithm is inferior to the
Newton method (and to BFGS giving the same results). But we can easily exploit  
the quadratic structure in Algorithm~\ref{A-Main}.\ref{BenchNewton}
by choosing the parameters properly and obtain the
exact solution after only $n$ steps too.

\subsection{Rosenbrock}
As next classical benchmark function we consider the Rosenbrock function
$f:\R^2\to\R$ given by
\begin{equation}\label{Rosenbrock}
 f(x,y)=(1-x)^2+100(y-x^2)^2\,.\tag{Ros}
\end{equation}
Here a steepest descent method leads into a canyon and then follows it, 
which is very time consuming. 
But the Newton method and Newton
based methods like BFGS show a very good performance for this
function. The minimizer is $(1,1)$ with value $f(1,1)=0$.

In Alt \cite{Alt2002} two versions of BFGS, the conjugated gradient method
{\bf (CG)} and two versions of trust region methods {\bf (TRM)}
are applied to the Rosenbrock
function where always initial point $(-1.9,2.0)$ is used. We used the same 
initial point and applied both 
Algorithm~\ref{A-Main}.\ref{Benchbasic} (with $\ep_0:=1.5$) and  
Algorithm~\ref{A-Main}.\ref{BenchNewton} (where $A_k$ is the Hessian at
iteration point $x_k$). The following table compares the results from
\cite[Section~$4.10.4$]{Alt2002} with our simulations.

\bigskip

{\small\begin{tabular}{||l|c c c c c||}
\hline\hline
Algorithm &BFGS~$1$          &BFGS~$2$  &CG        &TRM~$1$
&TRM~$2$\\ 
\hline
Iterations & $24$             &          &$47$      &$31$               &$34$\\
Gradients &$25$              &$49$      &$69$      &$93$               &$33$\\
Value     &$1.85\cdot10^{-6}$&$<10^{-5}$&$<10^{-3}$&$3.18\cdot 10^{-9}$&$1.3\cdot10^{-15}$\\
\hline\hline
Algorithm&\multicolumn{3}{c}{Algorithm~\ref{A-Main}.\ref{Benchbasic}
}&\multicolumn{2}{c||}{Algorithm~\ref{A-Main}.\ref{BenchNewton}}\\ 
\hline
Iterations &$14$&$19$&$29$&\multicolumn{2}{c||}{$12$}\\
Gradients&$29$&$37$&$55$&\multicolumn{2}{c||}{$20$}\\
Value    &$1.74\cdot10^{-6}$&$2.25\cdot10^{-9}$&$1.26\cdot10^{-18}$&\multicolumn{2}{c||}{$0.0$}\\
\hline\hline
\end{tabular}}

\bigskip

Algorithm~\ref{A-Main} shows again very good results. In particular 
Algorithm~\ref{A-Main}.\ref{BenchNewton} is impressive by 
giving the exact solution $(1,1)$ after merely 
$12$ iterations (i.e. $12$ computations of the Hessian) and $20$
gradient computations.

\subsection{Hilbert function}
For the comparison of Algorithm~\ref{A-Main}.\ref{Benchbasic} with a Newton
like algorithm (BFGS), we consider the quadratic
Hilbert function $f:\R^n\to\R$ given by 
\begin{equation}
f(x)=x^T\cdot A\cdot x \qmq{with} A=\Big(\frac1{i+j-1}\Big)_{1\leq
  i,j\leq n}\in\R^{n\times n}\,.\tag{Hilbert} 
\end{equation}
It is strictly convex with the unique minimizer $x=0$. 
The so-called Hilbert matrix $A$ is very ill-conditioned. In practice,
however, rounding errors in the representation of $A$ seem to 
improve the actual condition (cf. \cite{SP75}). 
For quasi Newton methods, $f$ is perfect by its quadratic
form, but bad by its condition. 

In Spedicato \cite{SP75} the BFGS algorithm was
applied with starting point  
$$x_0=\Big(\,\frac 41,\frac 42,\frac 43,\ldots,\frac 4n\,\Big)$$ 
and dimensions 
$$n=10,\; n=40\:\:\:\mbox{ and}\:\:\:n=80\,.$$
We used Algorithm~\ref{A-Main}.\ref{Benchbasic} with
$\ep_0=\sqrt{n}$ and the Euclidean norm $\norm{\cdot}$.
The results are given in the following table. 

\bigskip

\begin{tabular}{||l|c|c|c||}
\hline\hline
\multicolumn{4}{||c||}{BFGS}\\
\hline\hline
$n$        & $10$              & $40$              & $80$              \\[2mm]
Iterations  & $13$              & $43$              & $83$              \\
Value       & $4\cdot10^{-11} $ & $5\cdot10^{-9} $ & $1\cdot10^{-10} $ \\
\hline\hline
\multicolumn{4}{||c||}{Algorithm~\ref{A-Main}.\ref{Benchbasic}}\\
\hline\hline
$n$        & $10     $           & $40$               & $80 $      \\[2mm]
Iterations  & $40  $              & $40  $             & $40 $              \\
Gradients   & $58  $              & $69  $             & $77$              \\
Value       & $7.0\cdot10^{-10} $ & $2.2\cdot10^{-10}$ & $3.8\cdot10^{-10} $
\\[2mm] 
Iterations & $100 $              & $100 $             & $100$              \\
Gradients  & $123 $              & $176 $             & $194$              \\
Value      & $4.6\cdot10^{-13} $ & $3.3\cdot10^{-14}$ & $3.0\cdot10^{-14} $ \\
\hline\hline
\end{tabular}

\bigskip

The (smooth) Newton method applied to the Hilbert function would
end up in the minimizer after just one step in the case of exact computation.
But this does not happen in practice and, thus, 
the results depend on the solving algorithm for linear equations. 
This aspect has been studied intensively for this problem and is 
beyond the scope of this paper. Therefore we do not apply 
Algorithm~\ref{A-Main}.\ref{BenchNewton} that would give again
the (smooth) exact Newton method. 

Since the BFGS algorithm is designed for quadratic like functions,
we did not expect these relatively good results of our computations
(depending, of course, on the initial radius $\ep_0$).
It seems that our algorithm can compensate the bad condition of $A$ 
by the treatment of neighborhoods that gives a certain robustness even for
smooth functions. Let us still mention that first simulations 
with Algorithm~\ref{A-Main} for contact
problems, where the stiffnes matrix of the strain energy is ill conditioned,
are quite promising.

\subsection{Nesterov's Chebyshev-Rosenbrock functions}

As in Lewis \& Overton \cite{LO08} 
we consider the function $\tilde{f}:\R^n\to\R$ suggested by Nesterov
\begin{equation*}
\tilde{f}(x):=\frac14(x_1-1)^2+\sum\limits_{i=1}^{n-1}(x_{i+1}-2x_i^2+1)^2
\end{equation*}
and the nonsmooth version
\begin{equation*}
\hat{f}(x):=\frac14(x_1-1)^2+\sum\limits_{i=1}^{n-1}|x_{i+1}-2x_i^2+1|\,
\end{equation*}
that are used there to test the BFGS algorithm. 
Obviously both functions have the unique minimizer 
$\overline{x}=(1,\ldots,1)$ with $\ti f(\ol{x})=\hat f(\ol{x})=0$.
According to \cite{LO08} the BFGS algorithm typically generates 
iteration points that rapidly approach the highly oscillating manifold
\begin{equation*}
M:=\set{x\in\R^n}{x_{i+1}=2x_i^2-1,\:\:i=1,\ldots,n-1}\,
\end{equation*}
and then roughly follows it to the minimizer $\ol{x}\in M$. Therefore it is
reasonable for tests to start near $M$ and, in particular, 
the starting point 
$$\hat{x}:=(-1,1,1,\ldots,1)\in M$$
is usually considered. More precisely, BFGS always uses $\hat x$ in the smooth
case and takes random starting points for the nonsmooth $\hat f$. 
We used $\hat x$ in the nonsmooth case and a small
perturbation of it in the smooth case. 

\bsk

For the $i$-th Chebyshev polynomial $T_i$ one has
$T_2\circ T_i=T_{2i}$ and $T_2(s)=2s^2-1$. Therefore
we get for any $x=(x_1,x_2,\ldots,x_n)\in M$ with $x_1\in [-1,\,1]$ that  
$$  
x_{i+1}=T_{2^i}(x_1)=\cos\big(2^i\arccos(x_1)\big)  
$$ 
where the representation by trigonometric functions can be found in
\cite{Cheney}. Thus 
$$M\cap[-1,\,1]^n=\set{\big(\,x_1, \cos(2^1\arccos(x_1)),\ldots,
\cos(2^{n-1}\arccos(x_1))\,\big)}{x_1 \in[-1,\,1]}\,.$$
Notice that the $i$-th Chebyshev polynomial $T_i$ oscillates $2^i-1$ times
between values $-1$ and $1$ on $[-1,1]$, 
i.e. it reaches $2^i-1$ times both values $-1$ and $1$. 
Normally all algorithms (ours and others)  
do not follow the entire manifold $M$, but skip some of the oscillations. The
amount of skipped oscillations essentially depends not only on the choice of
parameters (like e.g. $\ep_0$) but also on the used computer. 
Naturally this
stochastic noise increases with the dimension of $X$ and the related 
strong increase of oscillations. Therefore this benchmark problem is not
recommended to compare the speed of algorithms (cf. \cite{LO08}).

\subsubsection{Smooth version}

First we consider the smooth function $\ti f$ where we compare results 
of the BFGS algorithm taken form \cite{LO08} with our computations. 
For Algorithm~\ref{A-Main}.\ref{Benchbasic} we use the
Euclidean norm, $\ep_0:=0.5$, and the special choice 
$\ThN(x)=0.001\,\frac{x}{\ep_0}$. Algorithm~\ref{A-Main}.\ref{BenchNewton} 
is applied with $\ep_0=0.5$ and $A_k$ being the Hessian
at iteration point $x_k$. Moreover we take the slightly perturbed initial point
$(-1.05,1,\ldots 1)$ for Algorithm~\ref{A-Main}.\ref{BenchNewton}, since it
finishes after just one step at the minimizer for initial point
$(-1,1,\ldots 1)$. The next table presents the results. 

\bigskip

\begin{tabular}{||l|c|c|c|c|c|c||}
\hline\hline
$\tilde{f}$ \phantom{$\dst\ti{\ti A}$} &\multicolumn{2}{|c|}{BFGS} &
\multicolumn{2}{|c|}{Algorithm~\ref{A-Main}.\ref{Benchbasic}}&
\multicolumn{2}{|c||}{Algorithm~\ref{A-Main}.\ref{BenchNewton}}\\ 
\hline\hline
Dimension & $n=8$ & $n=10$                 & $n=8$             & $n=10$
& $n=8$   & $n=10$\\[2mm] 
Iterations& $\approx 6,700$&$\approx50,000$&$16,683$           & $223,639$           & $4,109$ & $31,600$  \\
Gradients & &                              &$124,040$          & $1,773,929$         & $4,779$ & $37,305$ \\
Value     & $<10^{-15}$     & $<10^{-15} $ & $4\cdot 10^{-26}$ & $6.4\cdot 10^{-16}$ & $0.0$   & $9.9\cdot 10^{-16}$ \\
\hline\hline
\end{tabular}

\bigskip

Let us still mention that Algorithm~\ref{A-Main}.\ref{BenchNewton} with 
the quite large radius $\ep_0=15$ for $n=10$ needs merely  
$19$ iterations and $27$ gradients to get an iteration point with 
value less than $10^{-15}$. Here the iterations jump over a big area where
$M$ highly oscillates and, this way, avoid the most difficult part of following
$M$. In fact we had to adapt parameters carefully in order to reach that 
our algorithm follows the oscillations of $M$. Summarizing it turns out 
that our algorithm is robust enough to follow $M$, but
it can also skip an awful region by working with appropriate  
resultants on a sufficiently large ball. In any case it finds the minimizer,
while normally Algorithm~\ref{A-Main}.\ref{BenchNewton} is significantly
faster than Algorithm~\ref{A-Main}.\ref{Benchbasic}.
We do not know how long the BFGS algorithm typically followed the oscillations
of $M$ and how much was skipped.

\subsubsection{Nonsmooth version} 
Now we consider the nonsmooth function $\hat{f}$. 
In \cite{LO08} the BFGS algorithm is applied to $\hat{f}$ with random initial
point. For $n=3$ it produces iteration points with value less than $10^{-8}$,
but for $n=4$ it usually breaks down.   
Algorithm~\ref{A-Main}.\ref{Benchbasic} with $\ep_0=0.5$ and 
starting point  $(-1,1,\ldots,1)$ succeeds up to dimension $n=6$ and is 
overburdened from $n=7$. 

\bigskip       

\begin{tabular}{||l|c|c|c|c||}
\hline\hline
$\hat{f}$ \phantom{$\dst\ti{\ti A}$}
&\multicolumn{4}{|c||}{Algorithm~\ref{A-Main}.\ref{Benchbasic}}\\ 
\hline\hline
Dimension  & $n=3$& $n=4$     & $n=5$    & $n=6$\\[2mm]
Iterations & $4,365$ & $25,766$           & $219,886$         & $>1,500,000$\\
Gradients  & $13,691 $   & $106,714$          & $1,124,623$       & \\
Value      & $2.6\cdot 10^{-15} $   & $3.8\cdot10^{-10}$ & $1.0\cdot10^{-8}$ &$<10^{-9}$\\
\hline\hline
\end{tabular}

\bigskip

Notice that the value of $\hat f$ near the minimizer is merely the modulus
of small numbers instead of its square for $\ti f$. This tells us that the
approximation of the minimizers should be of comparable 
quality in both cases, though the values are larger in the nonsmooth case.

\subsubsection{Approximating a critical point which is not a minimizer}

Let us still take a short look at the nonsmooth variant $f:\R^n\to\R$
with 
\begin{equation*}
 f(x)=\frac14|x_1-1|+\sum\limits_{i=1}^{n-1}\big|x_{i+1}-2|x_i|+1\big|\,
\end{equation*}
having the same minimizer $\overline{x}=(1,\ldots,1)$ with value
$f(\ol{x})=0$.
The BFGS method, the gradient sampling algorithm, and Algorithm~\ref{A-Main}
try to follow the set 
\begin{equation*}
 N=\set{x\in\R^n}{x_{i+1}=2|x_i|-1\,:\;\;i=1,\ldots,n-1}\;,
\end{equation*}
which is, in contrast to $M$, not a differentiable manifold. 
There are also points on $N$ where $f$ is not differentiable. 
In particular for $n=2$, function $f$ is not differentiable at $(0,-1)$ and
one has $0\in\partial f(0,-1)$ though it is not a minimizer. 
The gradient sampling algorithm and the BFGS algorithm converge to this point
for many starting points according to \cite{LO08}. For
Algorithm~\ref{A-Main}.\ref{Benchbasic} we observed that it always converges
to either the critical point $(0,-1)$ (e.g. for initial point $(-1,1)$) or 
the global minimizer$(1,1)$.

Summarizing all cases, our algorithm appears to be quite robust in the presence
of high oscillations.

\subsection{Chebyshev approximation by exponential sums}

The gradient sampling algorithm (GS), which is also based on the concept of 
generalized gradients on sets, was widely tested in 
Burke, Lewis \& Overton \cite{CS-Verfahren-Burke}.
Let us take the Chebyshev approximation by exponential sums as
the most simple problem investigated in \cite{CS-Verfahren-Burke}, but which
keeps computations comparable. For a given 
function $u:[1,10]\to\R$ we are looking for minimizers 
$(a,b)\in\R^m\times\R^m$ of the function $\bar f:\R^n\to\R$ ($n=2m$) given by 
\begin{equation*}
\bar{f}(a,b):=\Big\|u(\cdot)-\sum\limits_{j=1}^ma_je^{-b_j(\cdot)} \Big\|_{\infty}
\end{equation*}
(where $\|\cdot\|_\infty$ denotes the supremum norm and $b_j(\cdot)$ the
product of $b_j$ with the argument).

As in \cite{CS-Verfahren-Burke} we first discretize the problem by fixing 
equidistant grid points $t_i$ in $[1,10]$ with 
$$t_i=1+9\tfrac i{N} \qmq{for} 0\leq i\leq N:=2000
$$ 
and, thus, we are looking for
minimizers $(a,b)\in\R^n$ of 
\begin{equation*}
f(a,b):=\max\limits_{0\leq i\leq N} 
\Big|u(t_i)-\sum\limits_{j=1}^ma_je^{-b_jt_i}\Big|\,.
\end{equation*}
Let us mention that a slightly modified problem was solved in
\cite{CS-Verfahren-Burke} to find the minimizer of $\bar f$. 
However it appeared that the minimization of $f$ is sufficient 
and we thus disregarded a further adaption. 

Now we fix $u(t)=\frac 1t$ on $[1,10]$ such that, with the smooth
functions   
\begin{equation*}
h_i(a,b):=\frac 1{t_i}-\sum\limits_{j=1}^ma_je^{-b_jt_i}\,,
\end{equation*}
we have to minimize 
\begin{equation}\label{fab}
f(a,b):=\max\limits_{0\leq i\leq N}\max\sset{h_i(a,b),\;-h_i(a,b)}
\end{equation}
on $\R^n$.

First we present the results from \cite{CS-Verfahren-Burke} obtained 
with the gradient 
sampling algorithm. Here the initial point $a=b=0$ is used and 
$2n$ gradients had to be computed in each 
iteration.\footnote{ We recall that the
  gradient sampling algorithm needs at least $n+1$ gradients in each
  iteration.} Due to lack of precision the algorithm didn't succeed for
$n>8$. 

\bigskip

\begin{tabular}{||c|c|c|c||}
\hline\hline
\multicolumn{4}{||c||}{Gradient sampling algorithm}\\
\hline\hline
$n=2m$ & $f$                   & Iterations & Gradients\\[2mm]
$2$ & $8.55641\cdot10^{-2}$ & $42 $      & $168$    \\
$4$ & $8.75226\cdot10^{-3}$ & $63 $      & $504$    \\
$6$ & $7.14507\cdot10^{-4}$ & $166$      & $1,992$   \\
$8$ & $5.58100\cdot10^{-5}$ & $282$      & $4,512$   \\
\hline\hline
\end{tabular}

\bigskip

For the application of Algorithm~\ref{A-Main}.\ref{Benchbasic} we start with
some observation. If we interpret the graph of $f$ as the surface of some
mountains, the initial point $a=b=0$ is located at some ridge where the path of
steepest descent follows that ridge until a saddle point. This means for the
algorithm that all iteration points $(a_k,b_k)$ and the corresponding gradients 
$f_k'=(a_k',b_k')$ have the form 
$$a_{k,j}=a_{k,1},\;b_{k,j}=b_{k,1},\;a'_{k,j}=a'_{k,1},\;b'_{k,j}=b'_{k,1}
\qmz{for all} j=1,\dots,m
$$
(which can be shown analytically by an easy induction argument)
and the sequence $(a_k,b_k)$
converges to a saddle point $(\tilde{a},\tilde{b})$ with  
\begin{equation*}
\tilde{a} = \frac2n(\hat{a},\,\hat{a},\,\ldots\,,\,\hat{a})\in\R^m\,, \quad
\tilde{b} = \frac2n(\hat{b},\,\hat{b},\,\ldots\,,\,\hat{b})\in\R^m
\end{equation*}
(e.g. $(\hat{a},\hat{b})\approx(1.43,\,0.45)\in\R^2$ for $n=2$).
It is remarkable that our algorithm also in practice precisely follows that
ridge and ends up in the saddle point. This however means that we need some
slight perturbation to find the minimizer. Therefore we consider both a
perturbed initial point 
$$a=-0.001\cdot\big(0^2,\,2^2,\ldots,(n-2)^2\big) \qmq{and} 
b=0.001\cdot\big(1^2,3^2,\ldots, (n-1)^2\big)$$
and a perturbed function 
$$ \hat{f}(a,b)=f(a, B\,b) \qmq{for some regular}
B\in\R^{m\times m}\,.
$$ 
More precisely, we replace $h_i(a,b)$ in \reff{fab} with 
\begin{equation*}
\hat{h}_i(a,b):=\frac 1{t_i}-\sum\limits_{j=1}^m a_je^{-jb_jt_i}\,
\end{equation*}
(i.e. $b_j$ is substituted by $jb_j$) and minimize
\begin{equation}
\hat{f}(a,b):=\max\limits_{0\leq j\leq
  N}\max\sset{\hat{h}_j(a,b),\;-\hat{h}_j(a,b)}\,.
\end{equation}
Notice that this change of $f$ does not effect the minimal value and, thus,
our results can be compared with those of the gradient sampling algorithm.
The next table presents the results of Algorithm~\ref{A-Main}.\ref{Benchbasic}
with $\ep_0=5\sqrt{m}$ and the special choice $\NHI(x)=0.1\cdot x$. 
\bigskip

\hspace*{-6mm}
\begin{tabular}{||c|c|c|c||c|c|c||}
\hline\hline
\multicolumn{7}{||c||}{Algorithm~\ref{A-Main}.\ref{Benchbasic}}\\
\hline\hline
&\multicolumn{3}{c||}{perturbed initial point} & 
\multicolumn{3}{c||}{function $\hat{f}$}\\  
$n=2m$ & $f(x_k)$   & Iterations & Gradients&  $\hat{f}(x_k)$                   & Iterations & Gradients\\[2mm]
$2$ & $8.55641\cdot10^{-2}$ & $10$      & $21$         & $8.55641\cdot10^{-2}$ & $14 $      & $32 $    \\
$4$ & $8.75226\cdot10^{-3}$ & $44$      & $124$        & $8.75226\cdot10^{-3}$ & $36 $      & $118$    \\
$6$ & $7.14509\cdot10^{-4}$ & $95$      & $431$        & $7.14509\cdot10^{-4}$ & $90 $      & $442$   \\
$8$ & $5.57688\cdot10^{-5}$ & $406$     & $2,547$       & $5.57688\cdot10^{-5}$ & $381$      & $2,512 $   \\
$10$& $4.24248\cdot10^{-6} $& $2,757$    & $22,075$     & $4.24249\cdot10^{-6} $& $2,766$     & $24,066$    \\
$12$& $3.17295\cdot10^{-7}$& $14,276 $ & $ 140,700$ & $3.17285\cdot10^{-7}$& $117,329 $ & $ 2,529,382$    \\
$14$&$8.56\cdot10^{-7}$& $17,961$&$180,065 $& $3.17570\cdot10^{-7}$ & $6,114 $   & $62,298 $  \\
\hline\hline
\end{tabular}

\bigskip

The deliberate choice of gradients might be some reason that we need 
much less gradients than the gradient sampling method.
We can also handle higher dimensions and reach 
our computational limit due to rounding errors for $n=14$.

\subsection{ Nonlinear regression}

We finally consider a nonlinear regression problem that can be found in Alt
\cite[Section~$2.3.1$]{Alt2002}. 
More precisely, 
we want to minimize $f:\R^3\to \R$ given by 
$$f(x)=\sum\limits_{i=1}^{10}\left(x_1e^{i\cdot x_2}+x_3-\eta_i\right)^2$$
for values

\bsk

\qquad\qquad\begin{tabular}{|c| cccccccccc|}
\hline
$i$    & $1$   & $2 $  & $3$   & $4$    & $5$    & $6$    &  $7$  & $8$   & $9$   & $10$  \\
$\eta_i$ & $1.0$ & $1.1 $& $1.2$ & $1.35$ & $1.55$ & $1.75$ & $2.5$ & $3.0$ & $3.7$ & $4.5$ \\
\hline
\end{tabular}

\bsk

It turns out that the steepest descent method,
the Nelder-Mead method, and the BFGS method are only partially able to solve
this problem for the initial points $x_0=(0,0,0)$ and $x_0=(1,1,1)$. 
More precisely we get from \cite{Alt2002} that: 
\bgl
\item[--]
The steepest descent method stops for both initial points at some
useless point (i.e. it is far from being a minimizer). 
\item[--]
For initial point $x_0=(1,1,1)$, the BFGS algorithm and the
Nelder-Mead method also stop at some useless point. With 
initial point $x_0=(0,0,0)$, the BFGS algorithm reaches 
a point with value $0.0861942$ after $21$
iterations and the Nelder-Mead method is also successful 
with a slightly larger value. 
\el
We successfully applied Algorithm~\ref{A-Main}.\ref{Benchbasic}
and Algorithm~\ref{A-Main}.\ref{BenchNewton} 
with $A_k$ being the Hessian at iteration point $x_k$ and with $\ep_0=0.5$ 
to both initial points.  

\bigskip

\begin{tabular}{||l| cc|cc||}
\hline\hline
Algorithm & \multicolumn{2}{|c}{Algorithm~\ref{A-Main}.\ref{Benchbasic}} &
\multicolumn{2}{|c||}{Algorithm~\ref{A-Main}.\ref{BenchNewton}}\\
\hline
Initial point &$(0,0,0)$  &$(1,1,1)$ &$(0,0,0)$ &$(1,1,1)$\\
Iterations    &$56$       &$42$      &$70$      &$49$\\
Gradients     &$130$      &$102$     &$194$     &$137$\\
Value
&\multicolumn{2}{|c|}{$0.0861942$}&\multicolumn{2}{|c||}{$0.0861942$}\\ 
\hline\hline
\end{tabular}

\bigskip

We also tested Algorithm~\ref{A-Main} with other starting points and 
always obtained the same minimal value of $f$  
and the same minimizer which is approximately 
$$(0.270,\,0.269,\,0.592)\,.$$

Though Algorithm~\ref{A-Main} is slower than the BFGS algorithm for the
first initial point, our algorithm seems to be much more robust in the sense
that it always finds the minimal solution 
(notice that the actual minimizer is unknown).  
It is interesting to see that even in this smooth case
Algorithm~\ref{A-Main}.\ref{Benchbasic} shows better convergence than 
Algorithm~\ref{A-Main}.\ref{BenchNewton} with variable norms.

\subsection{Summary}

Summarizing we can say that Algorithm~\ref{A-Main}
found the minimizer not only in all cases where the others succeeded, but also 
in cases where others failed. Moreover, for many problems our new algorithm 
needed the least iterations and gradient computations, i.e. it was the fastest
in this sense. In particular in higher dimensions the potential of 
Algorithm~\ref{A-Main} became clearly visible. 

It turns out that Algorithm~\ref{A-Main} is certainly an 
alternative for the solution of nonlinear and especially nonsmooth
minimization problems. Stability and robustness of  the algorithm 
are very convincing. It not only avoids typical oscillations of classical
smooth schemes, it also follows highly oscillating descent paths of 
very nasty functions and it can precisely trace some mountain ridge. 

Let us mention that we did not yet exploit the full potential of the
algorithm, since e.g. the choices of $G$, $T_1$, $T_2$ and the step size
control based on $\ep_{k,i}$ are not yet optimized by a systematic 
investigation.\footnote{For
many of the considered benchmark problems slower decreasing functions 
like $\NHI(x)=K\sqrt{x}$ instead of $\NHI(x)=K x$ 
lead to better results for large dimensions. 
Using $\USg(x,y)>y$ we also got better results in many cases.}
Also efficient stopping  
criteria are not yet considered. Nevertheless the achieved results 
are a promising basis for the treatment of 
relevant variational problems. Simulations for the highly degenerate
eigenvalue problem of the 1-Laplacian will be presented in an upcoming paper.

\end{document}